\documentclass{siamart0516}
\usepackage{tikz}
\usepackage{amssymb}
\usepackage{color,soul}

\def\GH{G_{\textup{H}}}
\def\RH{R_{\textup{H}}}
\def\GN{G}
\def\RN{R}

\input{defs.sty}

\usepackage{lipsum}
\usepackage{amsfonts}
\usepackage{graphicx}
\usepackage{epstopdf}
\usepackage{algorithmic}
\ifpdf
  \DeclareGraphicsExtensions{.eps,.pdf,.png,.jpg}
\else
  \DeclareGraphicsExtensions{.eps}
\fi

\newcommand{\TheTitle}{Asymptotic and numerical analysis of a stochastic PDE model of volume transmission}
\newcommand{\ShortTitle}{Asymptotic/numeric analysis of volume transmission}
\newcommand{\TheAuthors}{Sean D. Lawley and Varun Shankar}

\headers{\ShortTitle}{\TheAuthors}

\title{{\TheTitle}\thanks{Submitted to the editors \today.
\funding{SDL was supported by the National Science Foundation (DMS-1944574, DMS-1814832, DMS-1148230). VS was also supported by the National Science Foundation (DMS-1521748 and CISE-CCF 1714844).}}}

\author{Sean D. Lawley\thanks{Department of Mathematics, University of Utah, Salt Lake City, UT 84112 USA (\email{lawley@math.utah.edu}).}
\and{Varun Shankar\thanks{School of Computing, University of Utah, Salt Lake City, UT 84112 USA (\email{shankar@cs.utah.edu}).}}}

\usepackage{amsopn}

\ifpdf
\hypersetup{
  pdftitle={\TheTitle},
  pdfauthor={\TheAuthors}
}
\fi

\begin{document}

\maketitle

\begin{abstract}
Volume transmission is an important neural communication pathway in which neurons in one brain region influence the neurotransmitter concentration in the extracellular space of a distant brain region. In this paper, we apply asymptotic analysis to a stochastic partial differential equation model of volume transmission to calculate the neurotransmitter concentration in the extracellular space. Our model involves the diffusion equation in a three-dimensional domain with interior holes that randomly switch between being either sources or sinks. These holes model nerve varicosities that alternate between releasing and absorbing neurotransmitter, according to when they fire action potentials. In the case that the holes are small, we compute analytically the first two nonzero terms in an asymptotic expansion of the average neurotransmitter concentration. The first term shows that the concentration is spatially constant to leading order and that this constant is independent of many details in the problem. Specifically, this constant first term is independent of the number and location of nerve varicosities, neural firing correlations, and the size and geometry of the extracellular space. The second term shows how these factors affect the concentration at second order. Interestingly, the second term is also spatially constant under some mild assumptions. We verify our asymptotic results by high-order numerical simulation using radial basis function-generated finite differences.
\end{abstract}


\section{Introduction}
\label{sec:intro}


Synaptic transmission (ST) is a fundamental neural communication mechanism. In ST, a neuron fires an action potential that travels down its axon to a synapse that is adjacent to a second neuron. Upon arrival at the synapse, this electrical signal is translated into a biochemical response that results in neurotransmitters being released into the synaptic cleft. The neurotransmitters then rapidly diffuse across the synaptic cleft and bind to receptors on the second neuron, which then makes the second neuron more or less likely to fire, depending on whether the neurotransmitters are excitatory or inhibitory. Notice that ST is a ``one-to-one'' mode of communication, as it occurs through an essentially private channel between two neighboring neurons \cite{agnati2014}.

Complementing ST, there is an important ``many-to-many'' neural communication pathway called volume transmission (VT) \cite{fuxe2015volume, fuxe2016volume, fuxe10}. In VT, a collection of neurons projects to a distant brain volume (a nucleus or part of a nucleus) and increases the neurotransmitter concentration in the extracellular space in the distant volume by firing action potentials \cite{fuxe10}. The effect of this process is the modulation of ST in the distant volume, and for this reason VT is sometimes referred to as neuromodulation. Two important examples of VT are the serotonin projection from the dorsal raphe nucleus to the striatum \cite{blandina89, bonhomme95} and the dopamine projection from the substantia nigra to the striatum \cite{brooks00} (see \cite{fuxe10} for more examples). VT is thought to play an important role in maintaining the sleep-wakefulness cycle, motor control, and treating Parkinson's disease and depression \cite{fuxe10}. Furthermore, very recent experiments suggest that VT may be a common neurobiological mechanism for the control of feeding and other behaviors \cite{noble2018}.

In spite of the biological and medical significance of VT, mathematical modeling has only recently been used to study this process \cite{best2017,lawley18prob,lawley16neural}. In contrast, mathematical modeling has been used to study ST for decades. Indeed, the model of Hodgkin and Huxley transformed the field of neuroscience (Hodgkin and Huxley won the 1963 Nobel Prize in Physiology or Medicine) and continues to influence current research \cite{catterall2012,hodgkin52}. Furthermore, modeling ST has led to innovations in dynamical systems theory \cite{ermentrout12} and still prompts new mathematical questions \cite{goldwyn2011,handy2018,lawley2018hh}.

In this paper, we use a stochastic partial differential equation (PDE) model of VT to investigate the average neurotransmitter concentration in the extracellular space. In particular, we are interested in how the average neurotransmitter concentration varies across the space and how it depends on (i) the number, arrangement, and shape of nerve varicosities, (ii) neural firing statistics, (iii) the amount of neurotransmitter released when a neuron fires, (iv) the neurotransmitter diffusivity, and (v) the size and geometry of the extracellular space. We answer these questions in the case that the nerve varicosities are small compared to the distance between them. We discuss the biological relevance of our model in the Discussion section.

\begin{figure}
\begin{center}
\includegraphics[width=.6\linewidth]{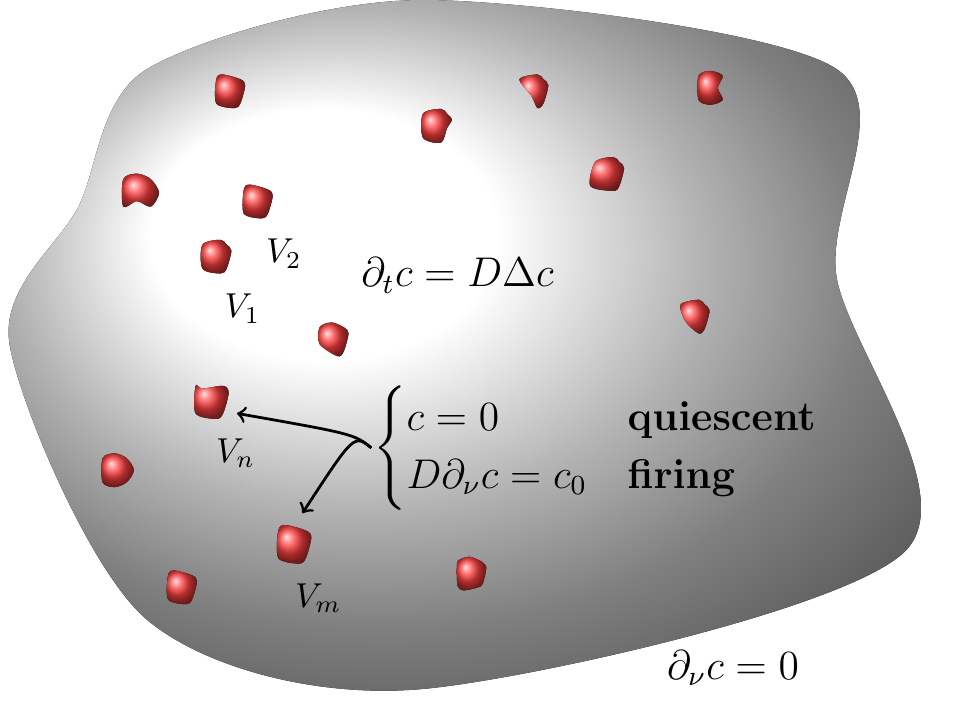}
\caption{Model schematic diagram. The three-dimensional grey region, $U\subset\mathbb{R}^3$, contains $N$ nerve varicosities, $\{V_n\}_{n=1}^{N}$, which are depicted by small red holes of arbitrary smooth shape. The concentration of neurotransmitter, $c(x,t)$, satisfies the diffusion equation in the extracellular space, $E=U\backslash\overline{\cup_{n=1}^NV_n}$, with boundary conditions at nerve varicosities that stochastically switch between absorption and flux into the space, corresponding to neurons firing action potentials. A no flux condition is imposed at the outer boundary.}\label{fig1}
\end{center}
\end{figure}

Mathematically, we use the diffusion equation in a bounded three-dimensional domain to model the neurotransmitter concentration in the extracellular space. This domain has a set of interior holes that represent nerve varicosities, see Figure~\ref{fig1}. Varicosities are bulbous enlargements on an axon that contain neurotransmitter \cite{goyal2013}. When a neuron fires, neurotransmitter is released from a nerve varicosity into the extracellular space, and when a neuron is not firing (quiescent), neurotransmitter is taken back up into the varicosity from the extracellular space. Hence, each of these interior holes has boundary conditions (either an inhomogeneous flux condition or an absorbing condition) that switch according to neural firing, which we assume follows a continuous-time Markov jump process. 

We first apply matched asymptotic analysis to a large system of deterministic PDEs that describes the steady-state mean of the stochastic PDE. Our asymptotic analysis yields the first two nonzero terms in an asymptotic expansion of the average neurotransmitter concentration, where the small parameter compares the size of varicosities to the distance between them. The first term in this expansion is constant in space, and this constant depends on very few details in the problem. Interestingly, if all the neurons have the same size, shape, and proportion of time firing, then the second term in the expansion is also constant in space. We verify our asymptotics by comparing to a simple example which can be solved analytically and by comparing to detailed numerical solutions for more complicated examples. 

We now comment on how the present study relates to prior work. A similar stochastic PDE model of VT in one space dimension was proposed and analyzed in \cite{lawley16neural}. This model was generalized to two and three space dimensions in \cite{lawley18prob}. In \cite{lawley18prob}, probabilistic methods were used to prove that the average neurotransmitter concentration is constant in space to leading order for small varicosities and to compute this constant. The asymptotic methods of the present study allow us to compute the next order term for the average neurotransmitter concentration and to relax the assumption in \cite{lawley18prob} that the varicosities are spherical. In further comparison to \cite{lawley18prob}, in the present study we develop numerical methods to study this stochastic PDE model.

More broadly, the present study employs and develops techniques from three disparate areas of mathematics. First, we use a recently developed theory that yields deterministic boundary value problems that describe statistics of parabolic PDEs with stochastically switching boundary conditions \cite{lawley15sima,pb1,lawley16bvp}. This theory has been used in diverse areas of biology, including insect respiration \cite{pb5}, intracellular reaction kinetics \cite{pb4}, and transport through gap junctions \cite{bressloff2016,pb6,pb13}. Second, our asymptotic analysis uses strong localized perturbation theory \cite{ward1993,ward1993log,ward2018}, which applies to perturbations of large magnitude but small spatial extent in elliptic PDE problems and reaction-diffusion systems. This theory has been used to study mean first passage times in the context of molecular and cellular biology \cite{cheviakov11,ward10,lindsay2017opt,pb3,lawley2018mfpt}, pattern formation in biological development \cite{ward09,tzou2017}, ecological dynamics \cite{kurella2015}, and diffusive signaling and chemoreception \cite{lindsay2017,bernoff2018,bernoff2018planar,lawley2018bp}. Finally, since the domain geometries considered in this article are irregular (in the sense of containing curved boundaries), we numerically discretize our spatial differential operators on scattered node sets using Radial Basis Function (RBF)-based finite difference (RBF-FD) formulas. RBF-FD formulas have been used in a variety of applications and are widely known to provide both numerical stability and high orders of convergence~\cite{Bayona2010,Davydov2011,Wright200699,FlyerNS,FlyerPHS,FlyerElliptic,ShankarJCP2017,SFJCP2018,FFBook,FF15}. 

The rest of the paper is organized as follows. We first briefly summarize our main results in the subsection below. In section~\ref{sec:dimensionless}, we give a precise and dimensionless formulation of our stochastic PDE model. In section~\ref{sec:mean}, we use asymptotic analysis to study the mean of this stochastic PDE. Section~\ref{sec:examples} applies our general asymptotic results to some examples. We describe our numerical methods and use them to verify our asymptotic results in section~\ref{sec:numerics}. We conclude by discussing the biological relevance of our model and describing future work.

\subsection{Main results}\label{main results}

Let $c(x,t)$ be the concentration of a neurotransmitter in the extracellular space $E\subset\R^{3}$ in some region of the brain. We assume that the neurotransmitter diffuses in the extracellular space,
\begin{alignat*}{2}
\frac{\partial c}{\partial t}&=D\Delta c,\quad &&x\in E\subset\R^{3},\;t>0.
\end{alignat*}
Furthermore, we assume that a set of neurons projects to this brain region, so that the domain $E$ contains $N\ge1$ holes in its interior which represent the corresponding nerve varicosities, see Figure~\ref{fig1}.

Each neuron fires action potentials and thus switches between a quiescent state and a firing state. We suppose that neuron firing is controlled by a continuous-time, reversible Markov jump process. We make no specific assumptions about firing correlations between neurons, and thus the neurons may fire synchronously, independently, or with some nontrivial correlations. When a neuron is not firing (quiescent), it absorbs neurotransmitter. When a neuron is firing, it releases neurotransmitter. Therefore, we impose the following stochastically switching boundary conditions at nerve varicosities,
\begin{alignat*}{2}
c(x,t)&=0,\quad &&x\in\partial V_{n},\quad\text{if neuron $n$ is quiescent at time }t,\\
D\partial_{\nu}c(x,t)&={\varphi}>0,\quad &&x\in\partial V_{n},\quad\text{if neuron $n$ is firing at time }t,
\end{alignat*}
where $\partial_{\nu}$ is the outward unit normal derivative and $\partial V_{n}$ is the boundary of the $n$-th nerve varicosity. Here, ${\varphi}>0$ is some constant flux of neurotransmitter into the extracellular space when the neuron is firing. We assume that the neurotransmitter does not leave the region through the outer boundary, and so we impose a no flux boundary condition there,
\begin{align*}
\partial_{\nu}c(x,t)=0,\quad\text{if }x\in\partial E\backslash\cup_{n=1}^{N}\partial V_{n}.
\end{align*}
Summarizing, our model consists of the diffusion equation on a three-dimensional domain with interior holes that randomly switch between being sources and sinks (this model was first considered in \cite{lawley18prob} with spherical varicosities).

Let $l_{0}$ denote the characteristic length scale describing the size of each nerve varicosity, and let $l_{1}>0$ denote the characteristic distance between varicosities. Assuming $\eps:=l_{0}/l_{1}\ll1$, we derive the following asymptotic expansion for the large time expected neurotransmitter concentration,
\begin{align}\label{main0}
\lim_{t\to\infty}\E[c(x,t)]
&=\frac{{\varphi}l_{1}}{D}\Big[\eps \theta^{(1)}
+\eps^{2}\Big(\theta^{(2)}+\Sigma(x)\Big)+\O(\eps^{3})\Big],\quad \eps\ll1.
\end{align}
The $\O(\eps)$ term is
\begin{align*}
\theta^{(1)}
&=\frac{\sum_{n=1}^{N}D_{n}\pi_{1}^{(n)}}{\sum_{n=1}^{N}C_{n}\pi_{0}^{(n)}},
\end{align*}
where $\pi_{1}^{(n)}\in(0,1)$ is the proportion of time that the $n$-th nerve varicosity is firing, $\pi_{0}^{(n)}=1-\pi_{1}^{(n)}$, and $C_{n}$ and $D_{n}$ are dimensionless constants determined by the shape of the $n$-th varicosity (see \eqref{ffar} and \eqref{Dnsa}). Notice that the $\O(\eps)$ term is spatially constant and independent of the following factors: (i) the positions of the nerve varicosities, (ii) the number of nerve varicosities (provided the proportion of varicosities of a certain shape and firing probability is fixed), (iii) the size and geometry of the extracellular space, (iv) firing correlations between the neurons, and (v) the rate of switching between firing and quiescent states. To see point (ii), notice that $\theta^{(1)}$ can be written as
\begin{align}\label{indepN}
\theta^{(1)}
=\frac{\mu_{D}}{\mu_{C}},
\end{align}
where $\mu_{D}$ and $\mu_{C}$ are the averages,
\begin{align*}
\mu_{D}
=\frac{1}{N}\sum_{n=1}^{N}D_{n}\pi_{1}^{(n)},\quad
\mu_{C}
=\frac{1}{N}\sum_{n=1}^{N}C_{n}\pi_{0}^{(n)}.
\end{align*}
Hence, if the number $N$ of nerve varicosities changes, but the proportion of varicosities of a certain shape and firing probability is fixed, then $\mu_{D}$ and $\mu_{C}$ are unchanged, and thus $\theta^{(1)}$ is unchanged.

The $\O(\eps^{2})$ term in \eqref{main0} shows how these factors affect the average neurotransmitter concentration. This term is the sum of a constant $\theta^{(2)}$ and function $\Sigma(x)$. For simplicity, we omit the formulas for $\theta^{(2)}$ and $\Sigma(x)$ from this section (we give the formula in section~\ref{sec:mean} below, see \eqref{theta2}-\eqref{sigmadef}). However, we note that if all the varicosities have the same firing probability and shape ($\pi_{1}^{(n)}=\pi_{1}^{(m)}$, $C_{n}=C_{m}$, and $D_{n}=D_{m}$ for all $n,m\in\{1,\dots,N\})$, then $\Sigma(x)\equiv0$, and thus this term is spatially constant.

\section{Problem setup}
\label{sec:dimensionless}

We now describe a dimensionless version of the model in section~\ref{main results} in further detail. Let $U\subset\R^{3}$ be open, connected, and bounded with a smooth boundary, $\partial U$. To represent the varicosities (or holes), let $\{V_{n}\}_{n=1}^{N}$ be $N$ open subsets of $U$, each with smooth boundary $\partial V_{n}$ and ``radius'' $\O(\eps)$. That is, we assume that there exists $N$ points $\{x_{n}\}_{n=1}^{N}$ in $U$ such that $V_{n}\to x_{n}$ uniformly as $\eps\to0$. We assume that the holes are well-separated in the sense that
\begin{align*}
|x_{n}-x_{m}|=\O(1)\quad\text{for }n\neq m,
\end{align*}
and
\begin{align*}
|x_{n}-x|=\O(1)\quad\text{for } n\in\{1,\dots,N\}\text{ and }x\in\partial U.
\end{align*}
We define the extracellular space, $E$, as the domain minus the holes, $$E:=U\backslash\overline{\cup_{n=1}^NV_n}.$$

To describe the state of each neuron (either firing or quiescent), let $\{J(t)\}_{t\ge0}$ be an irreducible continuous-time Markov jump process on a finite state space $\J\subseteq\{0,1\}^{N}$, where $J_{n}(t)\in\{0,1\}$ denotes its $n$-th component for $n\in\{1,\dots,N\}$. That is, $J(t)$ is the $N$-dimensional vector,
\begin{align*}
J(t)
=(J_{1}(t),\dots,J_{N}(t))\in\J\subseteq\{0,1\}^{N},
\end{align*}
with $J_{n}(t)\in\{0,1\}$ for $n\in\{1,\dots,N\}$. We emphasize that the state space $\J$ depends on the firing correlations between neurons. For example, if the neurons fire synchronously, then $J$ consists of only two elements, namely the vector of length $N$ with all $0$'s and the vector of length $N$ with all $1$'s,
\begin{align}\label{J2}
\J=\{(0,\dots,0),(1,\dots,1)\}\subsetneq\{0,1\}^{N}.
\end{align}
On the other hand, if the neurons fire independently, then $\J=\{0,1\}^{N}$.

Let $Q\in\R^{|\J|\times|\J|}$ denote the infinitesimal generator matrix of $J(t)$. Letting $Q_{\i,\j}\in\R$ denote the entry in the row corresponding to state $\i\in\J$ and the column corresponding to state $\j\in\J$, recall that $Q_{\i,\j}\ge0$ gives the rate that $J$ jumps from state $\i\in\J$ to state $\j\in\J$ for $\i\neq\j$, and the diagonal entries $Q_{\i,\i}$ are chosen so that $Q$ has zero row sums \cite{norris1998}.

Irreducibility ensures that $J(t)$ has a unique invariant distribution.
That is, there exists a unique $\rho\in\R^{|\J|}$ satisfying $Q^{T}\rho=0$ and $\sum_{\j\in\J}\rho_{\j}=1$, where $Q^{T}$ denotes the transpose of $Q$. We further assume that $Q$ is diagonalizable with all real eigenvalues. We note that if $J(t)$ is reversible (satisfies detailed balance), then these assumptions on $Q$ are satisfied. To see this, note that reversibility means 
\begin{align*}
\rho_{\i}Q_{\i,\j}
=\rho_{\j}Q_{\j,\i},\quad\text{for all }\i,\j\in\J.
\end{align*}
Hence, a straightforward calculation shows that $Q$ is self-adjoint with respect to the inner product
\begin{align*}
\langle\mu,\eta\rangle:=\sum_{\i\in\J} \rho_{\i}\mu_{\i}\eta_{\i},\quad \mu,\eta\in\R^{|\J|}.
\end{align*}
The spectral theorem thus ensures that $Q$ is diagonalizable with all real eigenvalues.

We note that assuming $Q$ is diagonalizable with all real eigenvalues is slightly restrictive. For example, if $N=3$ neurons always fire in a prescribed order (the first neuron fires, then the second, then the third, then the first, etc.), then a quick calculation shows that $Q$ has complex eigenvalues. However, $Q$ is diagonalizable with all real eigenvalues in many cases. For example, if the neurons fire synchronously or independently, then $J(t)$ is reversible (and thus $Q$ is diagonalizable with all real eigenvalues by the argument above). In fact, if the $N$ neurons can be partitioned into subgroups so that each subgroup fires synchronously and different subgroups fire independently, then $J(t)$ is reversible.

Assume that $\{c(x,t)\}_{t\ge0}$ is the $L^2(E)$-valued stochastic process that satisfies the diffusion equation in $E$ with a no flux condition at the outer boundary,
\begin{alignat}{2}
\frac{\partial c}{\partial t}&=\Delta c,\quad &&x\in E,\;t>0,\label{PDE}\\
\partial_{\nu}c&=0,\quad &&x\in\partial U,\label{outer}
\end{alignat}
and boundary conditions at each $V_n$ that switch according to the jump process $J_{n}(t)$,
\begin{alignat}{2}
c&=0,\quad &&x\in\partial V_n,\text{ if }J_n(t)=0,\label{0}\\
\partial_{\nu}c&=1,\quad &&x\in\partial V_n,\text{ if }J_n(t)=1.\label{1}
\end{alignat}
Since $c(x,t)$ is piecewise deterministic, it is straightforward to construct $c(x,t)$ by composing the deterministic solution operators corresponding to the  $|\J|$-many sets of boundary conditions according to the path of $\{J(t)\}_{t\ge0}$ (see \cite{lawley18ther} for a construction of a solution to a similar switching PDE).

This model is obtained from the model in section~\ref{main results} by defining dimensionless time, length, and concentration variables, $\frac{D}{l_{1}^{2}}t$, $\frac{1}{l_{1}}x$, and $l_{1}^{3}c$, for some length scale, $l_{1}>0$. In these variables, the dimensionless inhomogeneous flux boundary condition (when a neuron is firing) is $l_{1}^{4}{\varphi}/D$, but we have taken this to be unity in (\ref{1}) without loss of generality. In particular, if $c$ denotes the solution to \eqref{PDE}-\eqref{1}, then the solution to the problem with $l_{1}^{4}{\varphi}/D$ on the righthand side of \eqref{1} is simply $(l_{1}^{4}{\varphi}/D)c$ by linearity. Hence, the dimensional result in \eqref{main0} is obtained by solving \eqref{PDE}-\eqref{1}, then multiplying by $l_{1}^{4}{\varphi}/D$, and then multiplying by $l_{1}^{-3}$ to obtain a dimensional concentration.

\section{Asymptotic analysis of mean neurotransmitter}
\label{sec:mean}

\subsection{Decompose and diagonalize mean}

To analyze the large time mean of $c(x,t)$, we decompose the mean based on the state of the jump process $J$. Specifically, for each element $\j=(j_{1},\dots,j_{N})\in\J\subseteq\{0,1\}^{|\J|}$ of the state space of $J$, define
\begin{align}\label{vvj}
v_{\j}(x):=\lim_{t\to\infty}\E[c(x,t)1_{J(t)=\j}].
\end{align}
Here, $1_{J(t)=\j}$ denotes the indicator function,
\begin{align*}
1_{J(t)=\j}
=\begin{cases}
1 & \text{if }J(t)=\j,\\
0 &\text{if }J(t)\neq\j.
\end{cases}
\end{align*}
Putting the deterministic functions in \eqref{vvj} into a vector, $\v(x)=\{v_{\j}\}_{\j\in\J}\in\R^{|\J|}$, we have by Theorem~6 in \cite{lawley18prob} (see also \cite{lawley16bvp}) that
\begin{align}\label{vpde}
\Delta\v(x)+Q^{T}\v(x)=0,\quad x\in E,
\end{align}
where $E:=U\backslash\overline{\cup_{n=1}^NV_n}$ is $U$ minus the holes. Since there is a no flux boundary condition at the outer boundary regardless of the state of $J(t)$, we have that
\begin{align}\label{vnoflux}
\partial_{\nu}\v&=0,\quad x\in\partial U.
\end{align}
One advantage of the decomposition in (\ref{vvj}) is that $v_{\j}$ must satisfy the boundary condition at nerve varicosities corresponding to $J(t)=\j\in\J$,
\begin{align}\label{BCind}
\begin{split}
v_{\j}&=0,\quad x\in\partial V_{n}\quad\text{if }j_{n}=0,\\
\partial_{\nu}v_{\j}&=\rho_{\j},\quad x\in\partial V_{n}\quad\text{if }j_{n}=1,
\end{split}
\end{align}
where $\rho=\{\rho_{\j}\}_{\j\in\J}\in\R^{|\J|}$ is the unique invariant probability measure of $J(t)$. The boundary conditions \eqref{vnoflux}-\eqref{BCind} follow from  Theorem~6 in \cite{lawley18prob}. 

Writing (\ref{BCind}) in matrix notation, we have that
\begin{align}\label{dcb}
(I-B_{n})\v+B_{n}(\partial_{\nu}\v-\rho)&=0,\quad x\in\partial V_{n},
\end{align}
where $I\in\R^{|\J|\times|\J|}$ is the identity matrix and $B_{n}\in\{0,1\}^{|\J|\times|\J|}$ is the diagonal matrix with entries,
\begin{align}\label{amat}
(B_{n})_{\i,\j}=0\text{ if } \i\neq\j,
\quad\text{and}\quad
(B_{n})_{\j,\j}=j_{n}\in\{0,1\},\quad
\text{where }\j=(j_{1},\dots,j_{N}).
\end{align}
That is, $B_{n}$ is the $|\J|\times|\J|$ diagonal matrix whose diagonal entry corresponding to state $\j\in\J$ is 1 (respectively 0) if the $n$-th neuron is firing (respectively quiescent) when $J(t)=\j$. We emphasize that \eqref{dcb} is merely a compact way of writing the homogeneous Dirichlet conditions and inhomogeneous Neumann conditions that are given in \eqref{BCind}. Hence, while \eqref{dcb} may look like a Robin condition, it is in fact a vector of homogeneous Dirichlet conditions and inhomogeneous Neumann conditions due to the structure of $B_{n}$ given in \eqref{amat}.

Notice that $\v$ satisfies the coupled PDE \eqref{vpde} with decoupled boundary conditions \eqref{vnoflux} and \eqref{dcb}. We now decouple the PDE at the cost of coupling the boundary conditions in \eqref{dcb}.
 By assumption, we can diagonalize $Q^{T}$,
\begin{align}\label{evect}
P^{-1}Q^{T}P
=-\Lambda
=-\text{diag}(\lambda_{1},\lambda_{2},\dots,\lambda_{|\J|}),
\end{align}
where 
\begin{align*}
0=\lambda_{1}<\lambda_{2}\le\cdots\le\lambda_{|\J|},
\end{align*}
and the first column of $P\in\R^{|\J|\times|\J|}$ is the invariant measure $\rho\in\R^{|\J|}$, and the first row of $P^{-1}\in\R^{|\J|\times|\J|}$ is all 1's. Hence, setting
\begin{align*}
\u:=P^{-1}\v,
\end{align*}
we have that $\u$ satisfies
\begin{align}\label{uPDE}
\Delta\u-\Lambda\u=0,\quad x\in E,
\end{align}
with no flux boundary conditions at the outer boundary
\begin{align}\label{unoflux}
\partial_{\nu}\u=0,\quad x\in\partial U,
\end{align}
and boundary conditions at the interior holes given by
\begin{align}\label{uint}
(I-B_{n})P\u+B_{n}(P\partial_{\nu}\u-\rho)&=0,\quad x\in\partial V_{n}.
\end{align}

\subsection{Inner and outer expansions}

\begin{figure}
\begin{center}
\begin{tikzpicture}[scale=.45,every node/.style={minimum size=.8cm}]
\begin{scope}[shift={(-14,0)}]
\def\R{1} 
\def\scale{.075}
\def\pde{(2,-1.75)}
\def\pa{(-2.6,1.9)}
\def\pb{(3.1,1.7)}
\def\pc{(-2.1,-3.1)}
\def\pd{(-3.3,-1)}
\def\pe{(-3.2,1.1)}
\def\pf{(0,3)}
\def\pg{(2.5,2.3)}
\def\ph{(1.4,3.5)}
\def\pi{(3.74,.23)}
\def\pj{(-1.5,-.1)}
\def\ncolor{black}
\def\ecolor{black}
\def\esize{\large}

\fill [ball color=white] plot [smooth cycle, tension=.5] coordinates {(-4,-4) (-6,0) (-5,2) (-4,4) (0,5) (6,4) (6,1)  (7,-3) (0,-5)}; 

  \begin{scope}[shift={(-4.35,2)},scale=\scale]
\fill [ball color=\ncolor] circle (1);
  \end{scope}
  
    \begin{scope}[shift={(-4.65,-2)},scale=\scale]
\fill [ball color=\ncolor] circle (1);
  \end{scope}
  
      \begin{scope}[shift={(4,3.55)},scale=\scale]
\fill [ball color=\ncolor] circle (1);
  \end{scope}
  
      \begin{scope}[shift={(.5,-3.5)},scale=\scale]
\fill [ball color=\ncolor] circle (1);
  \end{scope}
  
        \begin{scope}[shift={(-3.7,-3.7)},scale=\scale]
\fill [ball color=\ncolor] circle (1);
  \end{scope}\  

	\begin{scope}[shift={(-3,3.5)},scale=\scale]
\fill [ball color=\ncolor] circle (1);
  \end{scope}

	\begin{scope}[shift={\pa},scale=\scale]
\fill [ball color=\ncolor] circle (1);
  \end{scope}
  	\begin{scope}[shift={\pc},scale=\scale]
\fill [ball color=\ncolor] circle (1);
  \end{scope}
  	\begin{scope}[shift={\pd},scale=\scale]
\fill [ball color=\ncolor] circle (1);
  \end{scope}
  	\begin{scope}[shift={\pe},scale=\scale]
\fill [ball color=\ncolor] circle (1);
  \end{scope}
    	\begin{scope}[shift={\pf},scale=\scale]
\fill [ball color=\ncolor] circle (1);
  \end{scope}
    	\begin{scope}[shift={\pg},scale=\scale]
\fill [ball color=\ncolor] circle (1);
  \end{scope}
    	\begin{scope}[shift={\ph},scale=\scale]
\fill [ball color=\ncolor] circle (1);
  \end{scope}
      	\begin{scope}[shift={\pi},scale=\scale]
\fill [ball color=\ncolor] circle (1);
  \end{scope}
      	\begin{scope}[shift={\pj},scale=\scale]
\fill [ball color=\ncolor] circle (1);
  \end{scope}

\node at \pde {{\esize\textcolor{\ecolor}{{\small $\Delta\u-\Lambda\u=0$}}}};
\node at (5,-4.5) {{\esize\textcolor{\ecolor}{{\small $\partial_{\nu} \u=0$}}}};

\node [below right  = -.2cm] at \pc {{{\small $x_m$}}};
\node [below right  = -.2cm] at \pd {{{\small $x_n$}}};
\node [below right  = -.2cm] at \pe {{{\small $x_1$}}};
\node [below right  = -.2cm] at \pa {{{\small $x_2$}}};

\node at (-6,4) {{{\small (a)}}};
\end{scope}


\begin{scope}[shift={(0,0)}]

\def\recsize{5.5}

\fill [color=gray!20] (-\recsize,-.85*\recsize) rectangle (\recsize,.85*\recsize);

\fill [ball color=red!35] plot [smooth cycle, tension=1] coordinates {(0,-2.5) (-2.5,0) (-2.1,2.1) (0,2.5) (2,2) (2.5,0) (1.5,-1)}; 

\draw[->] [black, thick] (0,0) -- (-2.1,2.1);

\node at (.5,1.5) {{{\small $|y_{n}|=\O(1)$}}};
\node at (0,4) {{\small\textcolor{black}{$\Delta_{y}\w_{n}-\eps^{2}\Lambda\w_{n}=0$}}};

\node at (0,-3) {{\scriptsize\textcolor{black}{$(I-B_{n})P\w_{n}+B_{n}(P\partial_{\nu}\w_{n}-\eps\rho)=0$}}};

\node at (-6,4) {{{\small (b)}}};

\end{scope}
\end{tikzpicture}
\caption{Construction of the matched asymptotic solution. (a) Outer solution $\u$ with a reflecting boundary condition on $\partial U$. From the perspective of the outer solution, the holes, $\{V_{n}\}_{n=1}^{N}$, have shrunk to the points, $\{x_{n}\}_{n=1}^{N}$. (b) Inner solution $\w_{n}$ in $\R^{3}\backslash U_{n}$ where $U_{n}:=\eps^{-1}V_{n}$.}\label{fig2}
\end{center}
\end{figure}
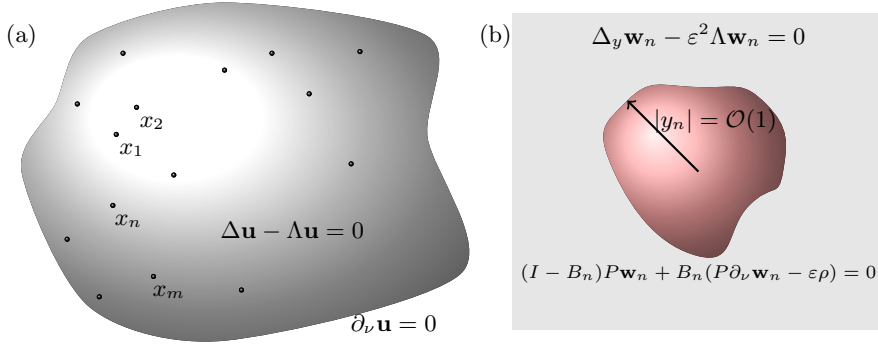

Using (\ref{evect}), observe that the first component of $\u$ (denoted by $(\u)_{1}$) is the large time mean
\begin{align}\label{sum}
(\u(x))_{1}
=\sum_{\j\in\J}v_{\j}(x)
=\lim_{t\to\infty}\E[c(x,t)].
\end{align}
Therefore, it remains to solve (\ref{uPDE}) subject to (\ref{unoflux}) and (\ref{uint}) in order to find the large time mean. We approximate the solution to this ${|\J|}$-dimensional boundary value problem in the case that the nerve varicosities are small, $\eps\ll1$. Specifically, we construct inner (or local) solutions that are valid in an $\O(\eps)$ neighborhood of each varicosity and then match these to an outer (or global) solution that is valid away from these neighborhoods, see Figure~\ref{fig2}. The asymptotic analysis that follows is formal and is related to the strong localized perturbation analysis pioneered in \cite{ward1993}.

The outer problem is
\begin{alignat*}{2}
\Delta\u-\Lambda\u&=0,\quad &&x\in U\backslash\cup_{n=1}^{N}\{x_{n}\},\\
\partial_{\nu}\u&=0,\quad &&x\in\partial U.
\end{alignat*}
Notice that from the perspective of the outer solution, the $N$ holes (varicosities), $\{V_{n}\}_{n=1}^{N}$, have shrunk to the points, $\{x_{n}\}_{n=1}^{N}$.
To construct an inner solution that is valid near the $n$-th nerve varicosity, we introduce the stretched coordinate
\begin{align*}
y:=\eps^{-1}(x-x_{n}),\quad n\in\{1,\dots,N\},
\end{align*}
and define $$\w_{n}(y;\eps):=\u(x_{n}+\eps y).$$ This inner solution satisfies
\begin{alignat}{2}
\Delta_{y}\w_{n}-\eps^{2}\Lambda\w_{n}&=0,&&\quad y\notin U_{n},\nonumber\\
\eps(I-B_{n})P\w_{n}+B_{n}(P\partial_{\nu}\w_{n}-\eps\rho)&=0,&&\quad y\in\partial U_{n},\label{wbcref}
\end{alignat}
where $U_{n}:=\eps^{-1}V_{n}$ is the magnified hole (note that $y\in\partial U_{n}$ if and only if $x\in\partial V_{n})$. Due to the structure of $B_{n}$ in \eqref{amat}, note that the factor of $\eps$ in the first term in the boundary condition in \eqref{wbcref} does not affect the solution. Hence, in the rest of the calculation we replace \eqref{wbcref} by 
\begin{align*}
(I-B_{n})P\w_{n}+B_{n}(P\partial_{\nu}\w_{n}-\eps\rho)&=0,\quad y\in\partial U_{n},
\end{align*}
since this does not affect the solution $\w_{n}$.

Next, we expand these inner and outer solutions in powers of $\eps$,
\begin{align*}
\u=\sum_{k=0}^{\infty}\eps^{k}\u^{(k)},
\quad
\w_{n}=\sum_{k=0}^{\infty}\eps^{k}\w_{n}^{(k)}.
\end{align*}
It follows that the $k$-th term in the outer solution expansion satisfies
\begin{alignat}{2}
\Delta\u^{(k)}-\Lambda\u^{(k)}&=0,\quad &&x\in U\backslash\cup_{n=1}^{N}\{x_{n}\},\label{vuk1}\\
\partial\u^{(k)}&=0,\quad &&x\in\partial U.\label{vuk2}
\end{alignat}
Similarly, the $k$-th term in the inner solution expansion satisfies
\begin{alignat}{2}
\Delta_{y}\w_{n}^{(k)}-\Lambda\w_{n}^{(k-2)}&=0,&&\quad y\notin U_{n},\label{vwk1}\\
(I-B_{n})P\w_{n}^{(k)}+B_{n}(P\partial_{\nu}\w_{n}^{(k)}-1_{k=1}\rho)&=0,&&\quad y\in\partial U_{n},\label{vwk2}
\end{alignat}
where $\w_{n}^{(-k)}\equiv0$ if $k>0$ and $1_{k=1}$ is 1 if $k=1$ and $0$ otherwise. The matching condition between the outer and inner solutions is that the near-field behavior of the outer solution as $x\to x_{n}$ must agree with the far-field behavior of the inner solution as $|y|\to\infty$. We write this formally as
\begin{align}
\u^{(0)}+\eps \u^{(1)}+\eps^{2}\u^{(2)}+\dots\sim \w_{n}^{(0)}+
\eps \w_{n}^{(1)}+\eps^{2}\w_{n}^{(2)}+\dots,\; x\to x_{n},\,|y|\to\infty\label{vmatch}.
\end{align}

\subsection{$\O(\eps^{0})$ term}\textcolor{white}{}
Consider the problem with no holes in the domain, $\eps=0$. Since $\lambda_{1} = 0$ (the first diagonal entry of $\Lambda$), the unperturbed outer solution, $\u^{(0)}(x)$, may be nontrivial in the first component, and we have the form
\begin{align}\label{u0vect}
\u^{(0)}(x)\equiv\Theta^{(0)}=(\theta^{(0)},0,\dots,0)\in\R^{|\J|},
\end{align}
for some constant $\theta^{(0)}\in\mathbb{R}$ to be determined.

With this information about $\u^{(0)}$, we can solve the inner problems for $\w_{n}^{(0)}$, $n\in\{1,\dots,N\}$. Let $f(y)$ be the real-valued solution to
\begin{align*}
\Delta_{y}  f&=0,
\quad
y\notin{U}_{n};
\qquad
f=1,
\quad
y\in\partial{U}_{n};
\qquad
f\to0,
\quad
|y|\to\infty.
\end{align*}
Then, the solution to the inner problem near the $n$-th hole is
\begin{align}\label{vw0}
\w_{n}^{(0)}(y)=\Theta^{(0)}-f(y)P^{-1}(I-B_{n})P\Theta^{(0)}.
\end{align}
To see this, observe that (\ref{vw0}) satisfies (\ref{vwk1}) since $f$ is harmonic. Further, (\ref{vw0}) satisfies (\ref{vwk2}) since by the structure of $B_{n}$ in \eqref{amat},
\begin{align*}
(I-B_{n})^{2}=I-B_{n},\quad
B_{n}(I-B_{n})=0,
\end{align*}
and $f=1$ on $\partial U_{n}$. Finally, (\ref{vw0}) satisfies (\ref{vmatch}) since $f\to0$ as $|y|\to\infty$.

We can now use (\ref{vw0}) to determine the leading order behavior of $\u^{(1)}$ as $x\to x_{n}$. The function $f(y)$ can be found in closed form only for a few simples shapes, but for a general shape it has the far-field asymptotic behavior
\begin{align}\label{ffar}
f\sim\frac{C_{n}}{|y|}+\frac{\Phi_{n}\cdot y}{|y|^{3}}+\dots,\quad
\text{as }|y|\to\infty,
\end{align}
where $C_{n}>0$ is called the electrostatic capacitance of ${U}_{n}$, and $\Phi_{n}\in\R^{3}$ is its dipole vector, and both are determined by $U_{n}$ \cite{jackson1975}. In general, $C_{n}$ must be calculated numerically, but it is known analytically for certain shapes, see Table~1 in \cite{cheviakov11}. In particular, if $U_{n}$ is a sphere, then $C_{n}$ is its radius. We will not need the value of $\Phi_{n}$.

Hence, (\ref{vw0}) and (\ref{ffar}) imply that $\w_{n}^{(0)}$ has the far-field behavior
\begin{align}\label{w0far}
\w_{n}^{(0)}\sim\Theta^{(0)}+\frac{\chi_{n}^{(1)}}{|y|},\quad\text{as }|y|\to\infty,
\end{align}
where
\begin{align}\label{eqf25}
\chi_{n}^{(1)}:=-C_{n}P^{-1}(I-B_{n})P\Theta^{(0)}.
\end{align}

Plugging (\ref{w0far}) into the matching condition (\ref{vmatch}) and using $y=\eps^{-1}|x-x_{n}|$, we have that
\begin{align*}
&\Theta^{(0)}+\eps \u^{(1)}+\dots
\sim 
\Theta^{(0)}+\frac{\eps \chi_{n}^{(1)}}{|x-x_{n}|}+\eps\w_{n}^{(1)}+\cdots.
\end{align*}
Hence, the outer solution $\u^{(1)}$ has the singular behavior as $x\to x_{n}$,
\begin{align*}
\u^{(1)}\sim \frac{\chi_{n}^{(1)}}{|x-x_{n}|},\quad x\to x_{n},\quad n\in\{1,\dots,N\}.
\end{align*}
We emphasize that $\u^{(1)}(x)$ is part of the expansion of the outer solution and is thus valid for $x$ outside an $\O(\eps)$ neighborhood of each of the holes. This singular behavior can be incorporated into the problem (\ref{vuk1}) for $\u^{(1)}$ in terms of the Dirac delta function as
\begin{align}\label{u0}
\Delta \u^{(1)}-\Lambda \u^{(1)}&=-4\pi\sum_{n=1}^{N}\delta(x-x_{n})\chi_{n}^{(1)},\quad x\in{U};\quad\partial_{\nu}\u^{(1)}=0,\quad x\in\partial U.
\end{align}
Integrating the first component of the $|\J|$-dimensional vector PDE in (\ref{u0}) and using the divergence theorem and the fact that the first row of $\Lambda$ is all zeros, we obtain that the sum over $n$ of the first components of $\chi_{n}^{(1)}$ must be zero,
\begin{align*}
\sum_{n=1}^{N}\big(\chi_{n}^{(1)}\big)_{1}=0.
\end{align*}
Now, using (\ref{amat}), (\ref{evect}), (\ref{u0vect}), and \eqref{eqf25}, we obtain
\begin{align*}
\big(\chi_{n}^{(1)}\big)_{1}
&=-C_{n}\sum_{\j\in\J}(1-j_{n})\rho_{\j}\theta^{(0)}.
\end{align*}
Therefore,
\begin{align}\label{vtheta0}
\theta^{(0)}=0,
\end{align}
and thus $\w_{n}^{(0)}=0$ by (\ref{vw0}).

\subsection{$\O(\eps^{1})$ term}\textcolor{white}{}
Since $\theta^{(0)} = 0$, we have that $\chi_{n}^{(1)} = 0$ for all $n\in\{1,\dots,N\}$ from \eqref{eqf25}. Thus, the solution to \eqref{u0} is the vector
\begin{align}\label{u1vect}
\u^{(1)}(x)\equiv\Theta^{(1)}=(\theta^{(1)},0,\dots,0)\in\R^{|\J|},
\end{align}
for some constant $\theta^{(1)}\in\mathbb{R}$ to be determined.

With this information about $\u^{(1)}$, we can solve the inner problem for $\w_{n}^{(1)}$, $n\in\{1,\dots,N\}$. Let $g(y)$ be the real-valued solution to
\begin{align}\label{g45}
\Delta_{y}  g&=0,
\quad
y\notin{U}_{n};
\qquad
\partial_{\nu}g=1,
\quad
y\in\partial{U}_{n};
\qquad
g\to0,
\quad
|y|\to\infty.
\end{align}
The solution to the inner problem near the $n$-th hole is
\begin{align}\label{vw1}
\w_{n}^{(1)}(y)=\Theta^{(1)}-f(y)P^{-1}(I-B_{n})P\Theta^{(1)}+g(y)P^{-1}B_{n}\rho.
\end{align}
To see this, observe that (\ref{vw1}) satisfies (\ref{vwk1}) since $f$ and $g$ are harmonic. Further, (\ref{vw1}) satisfies (\ref{vwk2}) since
\begin{align*}
(I-B_{n})^{2}=I-B_{n},\quad
B_{n}^{2}=B_{n},\quad
(I-B_{n})B_{n}=0=B_{n}(I-B_{n}),
\end{align*}
and $f=\partial_{\nu}g=1$ on $\partial U_{n}$. Finally, (\ref{vw1}) satisfies (\ref{vmatch}) since $f\to0$ and $g\to0$ as $|y|\to\infty$.

We can now use (\ref{vw1}) to determine the leading order behavior of $\u^{(2)}$ as $x\to x_{n}$. Using the known far-field behavior in \eqref{ffar} for solution to the analogous Dirichlet problem, we assume that the function $g(y)$ has the far-field asymptotic behavior
\begin{align}\label{gfar}
g\sim\frac{D_{n}}{|y|}+\frac{\Psi_{n}\cdot y}{|y|^{3}}+\dots,\qquad\text{as }|y|\to\infty,
\end{align}
where $D_{n}>0$ and $\Psi_{n}\in\R^{3}$ are determined by the shape of $U_{n}$. In fact, integrating the PDE in \eqref{g45} over a large sphere containing $U_{n}$ and using that $\partial_{\nu}g=1$ on $\partial U_{n}$ and \eqref{gfar}, we obtain
\begin{align}\label{Dnsa}
D_{n}
=\frac{|\partial U_{n}|}{4\pi},
\end{align}
where $|\partial U_{n}|$ is the surface area of the magnified hole, $U_{n}:=\eps^{-1}V_{n}$. We note that if $U_{n}$ is a sphere, then it thus follows that $D_{n}=C_{n}^{2}$, where $C_{n}$ is the radius of $U_{n}$. We will not need the value of $\Psi_{n}$.

Hence, $\w_{n}^{(1)}$ has the far-field behavior
\begin{align}\label{w1far}
\w_{n}^{(1)}\sim\Theta^{(1)}
+\frac{\chi_{n}^{(2)}}{|y|}
+\kappa_{n}\frac{\Phi_{n}\cdot y}{|y|^{3}}
+\mu_{n}\frac{\Psi_{n}\cdot y}{|y|^{3}}
,\qquad\text{as }|y|\to\infty,
\end{align}
where
\begin{align}\label{likef25}
\chi_{n}^{(2)}&:=D_{n}P^{-1}B_{n}\rho-C_{n}P^{-1}(I-B_{n})P\Theta^{(1)},\\
\kappa_{n}&:=-P^{-1}(I-B_{n})P\Theta^{(1)},\qquad \mu_{n}:=P^{-1}B_{n}\rho.\nonumber
\end{align}

Plugging (\ref{w1far}) into the matching condition (\ref{vmatch}) and using $y=\eps^{-1}|x-x_{n}|$, we have that
\begin{align}\label{refsugg}
\begin{split}
&\eps \Theta^{(1)}+\eps^{2}\u^{(2)}+\dots\\
&\quad\sim
\eps \Big(\Theta^{(1)}+\frac{\eps\chi_{n}^{(2)}}{|x-x_{n}|}+\kappa_{n}\frac{\eps^{2}\Phi_{n}\cdot (x-x_{n})}{|x-x_{n}|^{3}}
+\mu_{n}\frac{\eps^{2}\Psi_{n}\cdot (x-x_{n})}{|x-x_{n}|^{3}}\Big)+\eps^{2}\w_{n}^{(2)}+\cdots.
\end{split}
\end{align}
Hence, $\u^{(2)}$ has the singular behavior as $x\to x_{n}$,
\begin{align*}
\u^{(2)}\sim \frac{\chi_{n}^{(2)}}{|x-x_{n}|},\quad x\to x_{n},\quad n\in\{1,\dots,N\}.
\end{align*}
This singular behavior can be incorporated into the problem (\ref{vuk1}) for $\u^{(2)}$ in terms of the Dirac delta function as
\begin{align}\label{u2}
\Delta \u^{(2)}-\Lambda \u^{(2)}&=-4\pi\sum_{n=1}^{N}\delta(x-x_{n})\chi_{n}^{(2)},\quad x\in{U};\quad\partial_{\nu}\u^{(2)}=0,\quad x\in\partial U.
\end{align}
Integrating the first component of the $|\J|$-dimensional vector PDE in (\ref{u2}) and using the divergence theorem and the fact that the first row of $\Lambda$ is all zeros, we obtain $$\sum_{n=1}^{N}\big(\chi_{n}^{(2)}\big)_{1}=0.$$ Now, using (\ref{amat}), (\ref{evect}), (\ref{u1vect}), and \eqref{likef25}, we have that
\begin{align}\label{chi2}
\big(\chi_{n}^{(2)}\big)_{1}=D_{n}\pi_{1}^{(n)}-C_{n}\pi_{0}^{(n)}\theta^{(1)},
\end{align}
where we have defined the marginal probability that $J_{n}$ is 1 or 0,
\begin{align*}
\pi_{1}^{(n)}:=\sum_{\j\in\J}j_{n}\rho_{\j}=\lim_{t\to\infty}\P(J_{n}(t)=1),
\quad
\pi_{0}^{(n)}:=1-\pi_{1}^{(n)}.
\end{align*}
Therefore, 
\begin{align}\label{theta1}
\theta^{(1)}=\frac{\sum_{n=1}^{N}D_{n}\pi_{1}^{(n)}}{\sum_{n=1}^{N}C_{n}\pi_{0}^{(n)}}.
\end{align}
Equation~\eqref{theta1} generalizes a result proven in \cite{lawley18prob} using probabilistic methods for the special case that each hole is sphere.

\subsection{$\O(\eps^{2})$ term}
The solution to (\ref{u2}) can be written in terms of Green's functions. Let $\GN(x,\xi)$ denote the Neumann Green's function, which satisfies
\begin{align}\label{greens}
\begin{split}
\Delta \GN&=\frac{1}{|U|}-\delta(x-\xi),\quad x\in{U};
\qquad
\partial_{\nu} \GN= 0,\quad x\in\partial{U},\\
\GN(x,\xi)&=\frac{1}{4\pi|x-\xi|}+R(x,\xi);
\qquad
\int_{U}\GN(x,\xi)\,dx=0,
\end{split}
\end{align}
where $|U|$ denotes the volume of the domain $U$. For $\lambda>0$, let $\GH(x,\xi;\lambda)$ denote the modified Helmholtz Green's function, which satisfies
\begin{align}\label{hgreens}
\begin{split}
\Delta \GH(x,\xi;\lambda)-\lambda \GH(x,\xi;\lambda)&=-\delta(x-\xi),\quad x\in{U},\\
\partial_{\nu}\GH(x,\xi;\lambda)&=0,\quad x\in\partial{U},\\
\GH(x,\xi;\lambda)&= \frac{e^{-\sqrt{\lambda}|x-\xi|}}{4\pi|x-\xi|}+\RH(x,\xi;\lambda).
\end{split}
\end{align}
The functions $R$ and $\RH$ are called the regular parts of the Green's functions. In the special case that the domain $U$ is a sphere, $G$ and $\GH$ are known explicitly \cite{straube2009,cheviakov11}. We take advantage of these explicit formulas in see section~\ref{sec:examples} below. 

Putting these Green's functions into a matrix, let $\G(x,\xi)$ denote the $|\J|\times|\J|$ diagonal matrix,
\begin{align*}
\G(x,\xi)
=\text{diag}\Big(\GN(x,\xi),\GH(x,\xi;\lambda_{2}),\dots,\GH(x,\xi;\lambda_{|\J|})\Big)\in\R^{|\J|\times|\J|},
\end{align*}
where $(0,\lambda_{2},\dots,\lambda_{|\J|})$ are the diagonal entries of $\Lambda$. Hence, the solution to (\ref{u2}) is
\begin{align}\label{u2soln}
\u^{(2)}(x)=\Theta^{(2)}+4\pi\sum_{n=1}^{N}\G(x,x_{n})\chi_{n}^{(2)},
\end{align}
where $\Theta^{(2)}\in\R^{|\J|}$ is a vector of all zeros, except for the first component which is some constant $\theta^{(2)}\in\R$ to be determined.

Next, we use (\ref{u2soln}) to expand $\u^{(2)}$ as $x\to x_{n}$ and obtain
\begin{align*}
\u^{(2)}
&\sim\frac{\chi_{n}^{(2)}}{|x-x_{n}|}+\Theta^{(2)}+\Gamma_{n},\quad x\to x_{n},\quad n\in\{1,\dots,N\},
\end{align*}
where
\begin{align}\label{gamma}
\begin{split}
\Gamma_{n}
&:=-\sqrt{\Lambda}\chi_{n}^{(2)}
+4\pi\Big(\mathbf{R}(x_{n},x_{n})\chi_{n}^{(2)}+\sum_{m=1,m\neq n}^{N}\G(x_{n},x_{m})\chi_{m}^{(2)}\Big),
\end{split}
\end{align}
where $\Lambda$ is the diagonal matrix in \eqref{evect} and where we have defined the $|\J|\times|\J|$ diagonal matrix,
\begin{align*}
\mathbf{R}(x_{n},x_{n}):=\text{diag}\Big(\RN(x_{n},x_{n}),\RH(x_{n},x_{n};\lambda_{2}),\dots,\RH(x_{n},x_{n};\lambda_{|\J|})\Big)\in\R^{|\J|\times|\J|},
\end{align*}
containing the regular parts of the Green's functions.

Having obtained this asymptotic behavior of $\u^{(2)}$, we now return to the matching condition (\ref{vmatch}) to determine the far-field behavior of the inner solution, $\w_{n}^{(2)}$, $n\in\{1,\dots,N\}$. In particular, we find that
\begin{align*}
\w_{n}^{(2)}\sim\Theta^{(2)}+\Gamma_{n},\quad|y|\to\infty.
\end{align*}
Hence, since $\w_{n}^{(2)}$ must also satisfy (\ref{vwk1})-(\ref{vwk2}), we have that
\begin{align*}
\w_{n}^{(2)}(y)=\Theta^{(2)}+\Gamma_{n}-f(y)P^{-1}(I-B_{n})P(\Theta^{(2)}+\Gamma_{n}).
\end{align*}
Therefore, we obtain the next term in the far-field behavior of $\w_{n}^{(2)}$,
\begin{align*}
\w_{n}^{(2)}\sim\Theta^{(2)}+\Gamma_{n}+\frac{\chi_{n}^{(3)}}{|y|},\quad|y|\to\infty,
\end{align*}
where
\begin{align}\label{chi3ref}
\chi_{n}^{(3)}=-C_{n}P^{-1}(I-B_{n})P(\Theta^{(2)}+\Gamma_{n}).
\end{align}
Plugging this into the matching condition (\ref{vmatch}) (similar to the procedure that led to \eqref{refsugg}), we find that $\u^{(3)}$ has the leading-order singular behavior as $x\to x_{n}$,
\begin{align*}
\u^{(3)}\sim\frac{\chi_{n}^{(3)}}{|x-x_{n}|}+\kappa_{n}\frac{\Phi_{n}\cdot (x-x_{n})}{|x-x_{n}|^{3}}
+\mu_{n}\frac{\Psi_{n}\cdot (x-x_{n})}{|x-x_{n}|^{3}},\quad x\to x_{n},\quad n\in\{1,\dots,N\}.
\end{align*}
This singular behavior can be incorporated into the problem (\ref{vuk1}) for $\u^{(3)}$ in terms of the Dirac delta function as
\begin{align}\label{u3}
\begin{split}
\Delta \u^{(3)}-\Lambda \u^{(3)}&=-4\pi\sum_{n=1}^{N}\delta(x-x_{n})\chi_{n}^{(3)}
+4\pi\sum_{n=1}^{N}\kappa_{n}\Phi_{n}\cdot\nabla \delta(x-x_{n})\\
&\quad+4\pi\sum_{n=1}^{N}\mu_{n}\Psi_{n}\cdot\nabla \delta(x-x_{n})
,\quad x\in{U},
\end{split}
\end{align}
with no flux boundary conditions, $\partial_{\nu}\u^{(3)}=0,\;x\in\partial U$. Integrating the first component of the $|\J|$-dimensional vector PDE in (\ref{u3}) and using the divergence theorem and the fact that the first row of $\Lambda$ is all zeros, we obtain $$\sum_{n=1}^{N}\big(\chi_{n}^{(3)}\big)_{1}=0.$$ Using (\ref{amat}), (\ref{evect}) and \eqref{chi3ref}, we solve this equation for $\theta^{(2)}$ and find
\begin{align}\label{theta2}
\theta^{(2)}=-\frac{\sum_{n=1}^{N}C_{n}\sum_{\i,\j\in\J}(1-i_{n})P_{\i,\j}(\Gamma_{n})_{\j}}{\sum_{n=1}^{N}C_{n}\pi_{0}^{(n)}}.
\end{align}

\subsection{Summary}
\label{sec:mean_summary}
We now put these calculations together and make some observations. In the small hole limit, $\eps\to0$, the large time expected neurotransmitter concentration is given asymptotically in the outer region $|x-x_{n}|\gg \O(\eps)$ for $n\in\{1,\dots,N\}$ by
\begin{align}\label{casymp}
\lim_{t\to\infty}\E[c(x,t)]
\sim
\eps\theta^{(1)}+\eps^{2}\Big(\theta^{(2)}+\Sigma(x)\Big)+\O(\eps^{3}),\quad \eps\to0,
\end{align}
where $\theta^{(1)}$ and $\theta^{(2)}$ are given by (\ref{theta1}) and (\ref{theta2}), 
\begin{align}\label{sigmadef}
\Sigma(x)
=4\pi\sum_{n=1}^{N}\big(\chi_{n}^{(2)}\big)_{1}\GN(x,x_{n})
\end{align}
$\big(\chi_{n}^{(2)}\big)_{1}$ is (\ref{chi2}), and $\GN(x,x_{n})$ is the Neumann Green's function (\ref{greens}).

From (\ref{casymp}), first observe that the large time expected neurotransmitter concentration is constant in space at order $\eps$. Next, observe from (\ref{theta1}) that this constant is independent of (i) the positions of the nerve varicosities, $\{x_{n}\}_{n=1}^{N}$, (ii) number $N$ of nerve varicosities (provided the proportion of varicosities of a certain shape and firing probability is fixed, see \eqref{indepN}), (iii) the size and geometry of the extracellular space, $E$, and (iv) firing correlations between the neurons. Elaborating on (iv), this means that the constant is unchanged if the neurons fire synchronously (perfectly correlated), independently, or with other correlations. Finally, (\ref{theta1}) shows that this constant depends on the proportion of time spent firing but not on the overall rate of switching between firing and quiescent states.

From (\ref{casymp}) and (\ref{chi2}), further observe that if all the holes have a common shape (up to rotations) and firing probability, then $\big(\chi_{n}^{(2)}\big)_{1}=0$. Therefore in this case, observe from (\ref{casymp}) that the large time expected neurotransmitter concentration is in fact constant in space  at order $\eps^{2}$. Furthermore, from (\ref{casymp}), (\ref{theta2}), and \eqref{gamma}, observe that the order $\eps^{2}$ term describes how the number, positions, firing rates, and correlations of the neurons and the geometry of the extracellular space affects the mean neurotransmitter concentration. 

\section{Examples}
\label{sec:examples}

In this section, we consider some specific examples of firing correlations and geometries in order to illustrate and verify the general results of section~\ref{sec:mean}.

\subsection{Synchronous firing}\label{sync}

First suppose that all $N$ neurons fire synchronously. That is, suppose $J(t)\in\J$ where the state space $\J$ consists of the two elements (i.e.\ $|\J|=2$) given in \eqref{J2}. Namely, $\J$ consists of the $N$-dimensional vector of all 0's corresponding to every neuron being queiscent and the $N$-dimensional vector of all 1's corresponding to every neuron firing. If we label these two states by $0$ and $1$, then $J(t)$ has transition rates
\begin{align}\label{arrows}
\text{(quiescent)}\quad0\Markov{\alpha}{\beta}1\quad\text{(firing)}.
\end{align}
In this case, \eqref{vpde}-\eqref{BCind} simplify to the pair of PDEs,
\begin{align}\label{pair}
\begin{split}
0
&=\Delta v_{0}-\beta v_{0}+\alpha v_{1},\\
0
&=\Delta v_{1}+\beta v_{0}-\alpha v_{1},
\end{split}
\end{align}
with boundary conditions
\begin{align}\label{pairbc}
\begin{split}
\partial_{\nu}v_{0}
&=0,\quad\partial_{\nu}v_{1}=0,\quad\quad\,\,  x\in \partial U,\\
v_{0}
&=0,\quad\partial_{\nu}v_{1}=\tfrac{\beta}{\alpha+\beta},\quad x\in\partial V_{n},\quad n\in\{1,\dots,N\}.
\end{split}
\end{align}
It follows from \eqref{arrows} that
\begin{align}\label{matrices1}
B_{n}
=\begin{pmatrix}
0 & 0\\
0 & 1
\end{pmatrix},\;
P
=\begin{pmatrix}
\frac{\alpha}{\alpha+\beta} & 1\\
\frac{\beta}{\alpha+\beta} & -1
\end{pmatrix},
\;
P^{-1}
=\begin{pmatrix}
1 & 1\\
\frac{\beta}{\alpha+\beta} & -\frac{\alpha}{\alpha+\beta}
\end{pmatrix},
\;
\Lambda
=\begin{pmatrix}
0 & 0\\
0 & \alpha+\beta
\end{pmatrix}.
\end{align}
If we further suppose that each $V_{n}$ is a sphere of radius $\eps>0$ so that $C_{n}=D_{n}=1$ for $n\in\{1,\dots,N\}$, then
\begin{align}\label{t1s}
\theta^{(1)}=\beta/\alpha,
\quad
\big(\chi_{n}^{(2)}\big)_{1}=0,
\quad
\big(\chi_{n}^{(2)}\big)_{2}=-\frac{\beta}{\alpha+\beta},
\quad
\big(\Gamma_{n}\big)_{1}=0,
\end{align}
and thus
\begin{align*}
\big(\Gamma_{n}\big)_{2}
&=\frac{\beta}{\alpha+\beta}\Big[\sqrt{\alpha+\beta}-4\pi\Big(\RH(x_{n},x_{n},\alpha+\beta)+\sum_{m=1,m\neq n}^{N}\GH(x_{n},x_{m};\alpha+\beta)\Big)\Big].
\end{align*}
It then follows from \eqref{theta2} that
\begin{align}\label{t2e}
\begin{split}
\theta^{(2)}
&=-\frac{\alpha+\beta}{\alpha}\frac{1}{N}\sum_{n=1}^{N}(\Gamma_{n})_{2}\\
&=-\frac{\beta}{\alpha}\Big[\sqrt{\alpha+\beta}-\frac{4\pi}{N}\sum_{n=1}^{N}\Big(\RH(x_{n},x_{n},\alpha+\beta)+\sum_{m=1,m\neq n}^{N}\GH(x_{n},x_{m};\alpha+\beta)\Big)\Big].
\end{split}
\end{align}


Notice that our expression for $\theta^{(2)}$ (and thus our expression for the large time mean neurotransmitter in \eqref{casymp}) involves the Neumann Green's function of the modified Helmholtz equation \eqref{hgreens}. Since this Green's function cannot be computed analytically for a general domain $U\subset\R^{3}$, we now suppose that the domain $U$ is a sphere of radius $R_{0}>0$. In this case, the modified Helmholtz Green's function satisfying \eqref{hgreens} is given explicitly by \cite{straube2009}
\begin{align}\label{gs}
\GH(x,x_{0};\lambda)
=\frac{e^{-\sqrt{\lambda}|x-x_{0}|}}{4\pi|x-x_{0}|}+\RH(x,x_{0};\lambda),
\end{align}
where
\begin{align}\label{rs}
\begin{split}
\RH(x,x_{0};\lambda)
&=\frac{1}{4\pi}\sum_{n=0}^{\infty}\frac{a_{n}}{\sqrt{rr_{0}}}p_{n}(\cos\gamma)I_{n+1/2}(\sqrt{\lambda} r)I_{n+1/2}(\sqrt{\lambda} r_{0}),\\
a_{n}
&=(2n+1)\frac{(\sqrt{\lambda} R_{0})K_{n+3/2}(\sqrt{\lambda} R_{0})-nK_{n+1/2}(\sqrt{\lambda} R_{0})}{(\sqrt{\lambda} R_{0})I_{n+3/2}(\sqrt{\lambda} R_{0})+nI_{n+1/2}(\sqrt{\lambda} R_{0})}.
\end{split}
\end{align}
Here, $I_{\nu}(x)$ (respectively $K_{\nu}(x)$) is the modified Bessel function of the first (respectively second) kind of order $\nu$, $p_{n}(\cos\gamma)$ is the Legendre polynomial of degree $n$ with argument $\cos\gamma = \cos \theta \cos \theta_{0} + \sin \theta \sin \theta_{0} \cos(\varphi - \varphi_{0})$, and $(r, \theta, \varphi)$ (respectively $(r_{0}, \theta_{0}, \varphi_{0})$) denotes the spherical coordinates of $x$ (respectively $x_{0}$).

Summarizing, we have obtained the formula
\begin{align}\label{efo}
\lim_{t\to\infty}\E[c(x,t)]
\sim\eps\theta^{(1)}+\eps^{2}\theta^{(2)}+\O(\eps^{3}),\quad\text{as }\eps\to0,
\end{align}
where $\theta^{(1)}$ and $\theta^{(2)}$ have the analytical formulas \eqref{t1s} and \eqref{t2e}, where the Green's functions are given in \eqref{gs}-\eqref{rs}. In section~\ref{sec:numerics} below, we show that \eqref{efo} agrees with an independent numerical approximation to \eqref{pair}-\eqref{pairbc}.

\subsection{Spherically symmetric}\label{shell}

We now consider a simple example, which is the only example in which we can find an analytical formula for the mean neurotransmitter. We then use this example to further check our asymptotic result in \eqref{efo}.

Suppose there is a single spherical hole of radius $\eps>0$ located in the center of a sphere of radius $R_{0}>\eps$. That is, $N=1$, with
\begin{align*}
V_{1}:=\{x\in\R^{3}:|x|<\eps\}\subset U:=\{x\in\R^{3}:|x|<R_{0}\}.
\end{align*}
Let $\alpha$ and $\beta$ be as in the previous example. Finding the mean neurotransmitter reduces to solving \eqref{pair} in the annulus $|x|\in(\eps,R_{0})$ with boundary conditions
\begin{align}\label{pairbc0}
\begin{split}
\partial_{\nu}v_{0}
&=0,\quad\partial_{\nu}v_{1}=0,\quad\quad\,\, |x|=R_{0},\\
v_{0}
&=0,\quad\partial_{\nu}v_{1}=\tfrac{\beta}{\alpha+\beta},\quad|x|=\eps.
\end{split}
\end{align}
It is straightforward to use radial symmetry to find the constant value of the large-time mean neurotransmitter. In particular, letting $r=|x|$, diagonalizing the problem as in \eqref{uPDE}-\eqref{uint}, using the matrices in \eqref{matrices1}, and using radial symmetry gives that
\begin{align*}
\begin{pmatrix}
u_{0}\\
u_{1}
\end{pmatrix}
=\u
=P^{-1}\v
=P^{-1}\begin{pmatrix}
v_{0}\\
v_{1}
\end{pmatrix}
=\begin{pmatrix}
v_{0}+v_{1}\\
\tfrac{\beta}{\alpha+\beta}v_{0}
-\tfrac{\alpha}{\alpha+\beta}v_{1}
\end{pmatrix}
\end{align*}
satisfies
\begin{align}
0&=\tfrac{2}{r}\tfrac{d }{d r}u_{0}+\tfrac{d^{2} }{d r^{2}}u_{0},\quad r\in(\eps,R_{0}),\label{u0ode}\\
0&=\tfrac{2}{r}\tfrac{d }{d r}u_{1}+\tfrac{d^{2} }{d r^{2}}u_{1}-(\alpha+\beta)u_{1},\quad r\in(\eps,R_{0}),\label{u1ode}
\end{align}
with boundary conditions
\begin{align}
\tfrac{d }{d r}u_{0}
&=0,\quad \tfrac{d }{d r}u_{1}=0,\quad\quad\,\, r=R_{0},\label{bc1ref}\\
\tfrac{\alpha}{\alpha+\beta}u_{0}+u_{1}
&=0,\quad\tfrac{\beta}{\alpha+\beta}\tfrac{d }{d r}u_{0}-\tfrac{d }{d r}u_{1}=\tfrac{-\beta}{\alpha+\beta},\quad r=\eps.\label{bc2ref}
\end{align}
It follows directly from \eqref{u0ode} and the first boundary condition in \eqref{bc1ref} that $u_{0}$ is spatially constant. We can then solve \eqref{u1ode}, the second boundary condition in \eqref{bc1ref}, and the second boundary condition in \eqref{bc2ref} to obtain $u_{1}$ (since $\tfrac{d }{d r}u_{0}=0$). From this value of the function $u_{1}$, we then obtain the constant value of $u_{0}$ from the first boundary condition in \eqref{bc2ref}, which yields
\begin{align}\label{explshell}
u_{0}(|x|)
=v_{0}(x)+v_{1}(x)
=\frac{\beta}{\alpha}\frac{1}{\sqrt{\alpha+\beta}}h\big(\eps\sqrt{\alpha+\beta},R_{0}\sqrt{\alpha+\beta}\big),\quad\text{for all } x\in[\eps,R_{0}],
\end{align}
where 
\begin{align*}
h(a,b)&:=
 \Big[1+\frac{1}{a}-\frac{2(1+b)}{1+b+e^{2(b-a)}(b-1)}\Big]^{-1}.\label{explicit30}
\end{align*}
With this explicit expression, one can obtain the Taylor series expansion for $\eps\ll1$,
\begin{align}\label{taylor0}
v_{0}(x)+v_{1}(x)
=\frac{\beta}{\alpha}\eps
+\frac{\beta}{\alpha R_{0}}\left[\frac{1- R_{0}  \sqrt{\alpha +\beta } \tanh \left( R_{0}  \sqrt{\alpha +\beta }\right)}{1-(R_{0}  \sqrt{\alpha +\beta })^{-1}\tanh \left( R_{0}  \sqrt{\alpha +\beta }\right) }\right]\eps^{2}+\O(\eps^{3}).
\end{align}

We now verify that our more general asymptotic result in \eqref{efo} reproduces \eqref{taylor0} in this special case. It is immediate that $\theta^{(1)}$ in \eqref{t1s} agrees with \eqref{taylor0}. Furthermore, in this special case we have that $\theta^{(2)}$ in \eqref{t2e} reduces to
\begin{align}\label{t2e0}
\begin{split}
\theta^{(2)}
&=-\frac{\alpha+\beta}{\alpha}(\Gamma_{1})_{2}
=-\frac{\beta}{\alpha}\Big[\sqrt{\alpha+\beta}-4\pi \RH(0,0;\alpha+\beta)\Big].
\end{split}
\end{align}
Now, it follows from \eqref{rs} that
\begin{align}\label{RHLIMIT}
\lim_{|x_{1}|\to0}\RH(x_{1},x_{1};\alpha+\beta)
=\frac{\sqrt{\alpha+\beta }+(\alpha+\beta) R_{0}}{2 \pi  \left(\sqrt{\alpha+\beta } R_{0}+e^{2 \sqrt{\alpha+\beta } R_{0}} \left(\sqrt{\alpha+\beta } R_{0}-1\right)+1\right)}.
\end{align}
To derive \eqref{RHLIMIT}, let $r=|x_{1}|$ and note that \eqref{rs} yields
\begin{align}\label{RHsame}
\RH(x_{1},x_{1};\alpha+\beta)
&=\frac{1}{4\pi}\sum_{n=0}^{\infty}\frac{a_{n}}{r}p_{n}(1)\big(I_{n+1/2}(\sqrt{\alpha+\beta} r)\big)^{2}.
\end{align}
Now, it is a property of the modified Bessel function that
\begin{align*}
\lim_{r\to0+}\frac{1}{r}\big(I_{n+1/2}(\sqrt{\alpha+\beta} r)\big)^{2}
=\begin{cases}
\frac{2}{\pi}\sqrt{\alpha+\beta} & \text{if }n=0,\\
0 & \text{if }n>0.
\end{cases}
\end{align*}
Hence, taking the limit $|x_{1}|=r\to0+$ in \eqref{RHsame} yields only the $n=0$ term,
\begin{align}\label{RHpre}
\lim_{|x_{1}|\to0}\RH(x_{1},x_{1};\alpha+\beta)
=\frac{1}{4\pi}a_{0}p_{0}(1)\frac{2}{\pi}\sqrt{\alpha+\beta}
=\frac{1}{2\pi^{2}}a_{0}\sqrt{\alpha+\beta}.
\end{align}
Plugging the value $a_{0}=K_{3/2}(\sqrt{\alpha+\beta}R_{0})/I_{3/2}(\sqrt{\alpha+\beta}R_{0})$ from \eqref{rs} into \eqref{RHpre} and simplifying yields \eqref{RHLIMIT}.

Plugging \eqref{RHLIMIT} into \eqref{t2e0} yields that
\begin{align*}
\theta^{(2)}
=\frac{\beta}{\alpha R_{0}}\left[\frac{1- R_{0}  \sqrt{\alpha +\beta } \tanh \left( R_{0}  \sqrt{\alpha +\beta }\right)}{1-(R_{0}  \sqrt{\alpha +\beta })^{-1}\tanh \left( R_{0}  \sqrt{\alpha +\beta }\right) }\right],
\end{align*}
which agrees with \eqref{taylor0}.

\subsection{Independent firing}

Suppose that each of the $N\ge1$ nerve varicosities, ${V}_{n}$, fires according to an independent Markov jump process $J_{n}(t)$ with transition rates $\alpha_n,\beta_n$,
\begin{align*}
\text{(quiescent)}\quad0\Markov{\alpha_n}{\beta_n}1\quad\text{(firing)},\quad n\in\{1,\dots,N\},
\end{align*}
and $J(t)=(J_1(t),\dots,J_N(t))\in\J=\{0,1\}^N$. In this case,
\begin{align*}
Q_{\i\j}=\sum_{n=1}^Nq^{(n)}_{i_n,j_n}\prod_{m\ne n}\delta_{i_m,j_m},\quad \j\in\J,
\end{align*}
where $\delta_{i,j}$ denotes the Kronecker delta and $q_{i,j}^{(n)}$ are the entries in the generator $Q^{(n)}$ of the Markov process $J_{n}(t)$ controlling the firing of the $n$-th nerve varicosity,
\begin{align*}
Q^{(n)}=\begin{pmatrix}
q^{(n)}_{0,0} & q^{(n)}_{0,1}\\
q^{(n)}_{1,0} & q^{(n)}_{1,1}
\end{pmatrix}
=
\begin{pmatrix}
-\beta_n & \beta_n\\
\alpha_n & -\alpha_n
\end{pmatrix}.
\end{align*}
Observe that $Q_{\i\j}=0$ if $\i$ and $\j$ differ in more than one component, and either $Q_{\i\j}=\alpha_{n}$ or $Q_{\i\j}=\beta_{n}$ if $\i$ and $\j$ differ in exactly the $n$-th position. Further, $Q_{\i\i}=\sum_{n=1}^{N}q_{i_{n},i_{n}}<0$ with $q_{i_{n},i_{n}}=-(1-i_{n})\beta_{n}-i_{n}\alpha_{n}$. Also,
\begin{align*}
\rho_{\j}=\prod_{n=1}^{N}\frac{j_{n}\beta_{n}+(1-j_{n})\alpha_{n}}{\alpha_n+\beta_n},\quad \j\in\J.
\end{align*}
Since the components of $J(t)$ are independent, the probability that the $n$-th component is 1 is $\frac{\alpha_{n}}{\alpha_{n}+\beta_{n}}$ and thus
\begin{align*}
\theta^{(1)}=\frac{\sum_{n=1}^{N}D_{n}\frac{\beta_{n}}{\alpha_{n}+\beta_{n}}}{\sum_{n=1}^{N}C_{n}\frac{\alpha_{n}}{\alpha_{n}+\beta_{n}}}.
\end{align*}
If $\alpha_{m}=\alpha_{n}$ and $\beta_{m}=\beta_{n}$ for all $m,n\in\{1,\dots,N\}$ then $\theta^{(1)}$ is the same as if the neurons fired synchronously (as noted above in section~\ref{sec:mean_summary}). 

To simplify the calculation of $\theta^{(2)}$, we assume that there are only two varicosities, $N=2$, and that each one is a sphere of radius $\eps>0$ so that $C_{n}=D_{n}=1$ for $n=1,2$. Further suppose that each varicosity has the same firing rate so that $\alpha_{n}=\alpha$ and $\beta_{n}=\beta$ for $n=1,2$. In this case, $\theta^{(1)}=\beta/\alpha$. Furthermore, it follows that
\begin{align}\label{matrices2}
\begin{split}
Q^{T}
&=\begin{pmatrix}
-2\beta & \alpha & \alpha & 0\\
\beta & -(\alpha+\beta) & 0 & \alpha\\
\beta & 0 & -(\alpha+\beta) & \alpha\\
0 & \beta & \beta & -2\alpha
\end{pmatrix},\\
P
&=\frac{1}{(\alpha+\beta)^{2}}
\begin{pmatrix}
\alpha^{2} & -\alpha\beta & 0 & \beta^{2}\\
\alpha\beta & (\alpha-\beta)\beta & -\beta^{2} & -\beta^{2}\\
\alpha\beta & 0 & \beta^{2} & -\beta^{2}\\
\beta^{2} & \beta^{2} & 0 & \beta^{2}
\end{pmatrix},\\
P^{-1}
&=\begin{pmatrix}
1 & 1 & 1 & 1\\
-2 & -1+\tfrac{\alpha}{\beta} & -1+\tfrac{\alpha}{\beta} & \tfrac{2\alpha}{\beta}\\
1-\tfrac{\alpha}{\beta} & -\tfrac{2\alpha}{\beta} & 1+\tfrac{\alpha^{2}}{\beta^{2}} & \tfrac{\alpha(\alpha-\beta)}{\beta^{2}}\\
1 & -\tfrac{\alpha}{\beta} & -\tfrac{\alpha}{\beta} & \tfrac{\alpha^{2}}{\beta^{2}}
\end{pmatrix},
\end{split}
\end{align}
and 
\begin{align*}
\Lambda=\text{diag}(0,\alpha+\beta,\alpha+\beta,2(\alpha+\beta)),\quad
B_{1}=\text{diag}(0,1,0,1),\quad 
B_{2}=\text{diag}(0,0,1,1).
\end{align*}
Using \eqref{theta2} and \eqref{matrices2}, a tedious but straightforward calculation yields
\begin{align}\label{t2ei}
\theta^{(2)}
=\frac{\beta}{\alpha}
\Big[2\pi\big(\RH(x_{1},x_{1};\alpha+\beta)+\RH(x_{2},x_{2};\alpha+\beta)\big)-\sqrt{
\alpha+\beta}\Big],
\end{align}
where $x_{1}$ and $x_{2}$ are the centers of the two varicosities.

We now compare the case of independent firing considered here to the case of synchronous firing in section~\ref{sync} above. Let $\theta_{\text{indep}}^{(2)}$ denote the value of $\theta^{(2)}$ in \eqref{t2ei} for independent firing and $\theta_{\text{sync}}^{(2)}$ denote the value of $\theta^{(2)}$ in \eqref{t2e} for synchronous firing when $N=2$. Comparing \eqref{t2e} to \eqref{t2ei}, we have that
\begin{align}\label{is}
\theta^{(2)}_{\text{sync}}
=\theta_{\text{indep}}^{(2)}+\frac{\beta}{\alpha}4\pi \GH(x_{1},x_{2};\alpha+\beta).
\end{align}
To understand \eqref{is}, observe that the Green's function $\GH(x,\xi;\lambda)$ for $\lambda>0$ is the steady-state concentration at position $x\in U$ of a diffusing chemical that degrades at rate $\lambda>0$, has a source at position $\xi\in U$, and reflects from the boundary $\partial U$. It follows that (i) $\GH(x,\xi;\lambda)>0$, (ii) $\GH(x,\xi;\lambda)$ grows as $x$ approaches $\xi$, and (iii) $\GH(x,\xi;\lambda)$ grows as $\lambda$ decreases (these three intuitive properties can be proved by analyzing \eqref{hgreens}).

Using these three properties of $G$, equation~\eqref{is} implies the following three conclusions. First, the neurotransmitter concentration for synchronous firing is strictly greater than the neurotransmitter concentration for independent firing. Second, this difference between synchronous and independent firing is more pronounced when the varicosities are close to each other. Third, this difference is also more pronounced when the overall rate of firing, $\alpha+\beta$, is slow. Moreover, while \eqref{is} concerns perfectly synchronous versus perfectly independent firing for $N=2$ varicosities, we conjecture that these three conclusions hold when comparing more correlated versus less correlated firing of $N\ge2$ varicosities.

\section{Numerical analysis of mean neurotransmitter}
\label{sec:numerics}

In this section, we verify our asymptotic results through numerical simulations. For these simulations, we consider the case of synchronous firing described in sections~\ref{sync}-\ref{shell} above for three different domain geometries. First, we consider the spherically symmetric case of section~\ref{shell}. Since the large-time mean neurotransmitter can be found explicitly in this case (see \eqref{explshell}), this case serves to validate our numerical methods. We then move on to more complicated geometries that cannot be solved analytically to verify our asymptotic result in \eqref{efo}. 
Before detailing these numerical experiments and their results, we first describe the discretization of the PDE \eqref{pair} and boundary conditions \eqref{pairbc} on these domains.

\subsection{Spatial discretization}
\label{sec:space_disc}

As mentioned in section~\ref{sec:intro}, we use the RBF-FD method for spatial discretization. More specifically, we use a specific variant of the RBF-FD method (developed by the second author) called the \emph{overlapped RBF-FD method}~\cite{ShankarJCP2017,SFJCP2018}. This variant of the RBF-FD method uses a judicious mix of centered and one-sided differences on scattered nodes, leading to fewer finite difference stencils than collocation points. This allows for improved computational efficiency over standard RBF-FD without sacrificing high-order accuracy or numerical stability.  We now briefly describe the overlapped RBF-FD method for approximating the Laplacian; the same discussion applies to any linear operator (including the Neumann boundary condition operator).

Let $X = \{{{x}}_k\}_{k=1}^N$ be a set of collocation nodes on some domain $\Omega \subset \mathbb{R}^3$ with boundary $\partial \Omega$. Define the stencil $P_k$ to be the set of nodes containing the nodes ${{x}}_{\calI^k_1}$, $\{{{x}}_{\calI^k_2},\hdots,{{x}}_{\calI^k_n}\}$, \emph{i.e.}, ${{x}}_{\calI^k_1}$ and its $n-1$ nearest neighbors. Here, $\calI^k_j$, $j=1,\ldots n$ are indices that map stencil nodes into the set $X$. Without loss of generality, we confine our discussion to the stencil $P_1$ and the Laplacian $\Delta$. First, given a parameter $\delta \in (0,1]$, we can define a ball associated with stencil $P_1$ with center ${{x}}_{\calI^1_1}$ and radius $r_1$ given by:
\begin{align}
r_1 = \left(1-\delta\right) \max\limits_{j} \|{{x}}_{\calI^1_1} - {{x}}_{\calI^1_j}\|,\quad j=1,\hdots,n.
\end{align}
This allows us to define the subset $\mathbb{B}_1 \subset P_1$ as:
\begin{align}
\mathbb{B}_1 = \{{{y}} \in P_1 : \|{{y}} - {{x}}_{\calI^1_1}\|_2 \leq r_1 \},
\end{align} 
\emph{i.e.}, $\mathbb{B}_1$ contains all the nodes in $P_1$ that lie within the ball of radius $r_1$ centered at ${{x}}_{\calI^1_1}$. Finally, let the global indices of the $p_1$ nodes in $\mathbb{B}_1$ be given by:
\begin{align}
R_1 = \{\calR^1_1, \calR^1_2, \hdots, \calR^1_{p_1}\}.
\label{eq:ball_inds}
\end{align}
The overlapped RBF-FD method now involves using all the nodes in $P_1$ (${{x}} \in P_1$) to compute RBF-FD weights for the nodes in $\mathbb{B}_1$ (${{y}} \in \mathbb{B}_1$). To do so, we first define the augmented local RBF interpolant on $P_1$:
\begin{align}
s_1({{x}},{{y}}) = \sum\limits_{j=1}^n w^1_j({{y}}) \|{{x}} - {{x}}_{\calI^1_j}\|^m + \sum\limits_{i=1}^{M} \lambda^1_i({{y}}) \psi^1_i({{x}}),
\label{eq:rbf_interp}
\end{align}
where $\|{{x}} - {{x}}_{\calI^1_j}\|^m$ is the polyharmonic spline (PHS) RBF of degree $m$ ($m$ is odd), and $\psi^1_i({{x}})$ are the $M$ monomials corresponding to the total degree polynomial of degree $\ell$ in $d$ dimensions. In order to compute the RBF-FD weights $w^1_j({{y}})$ for the Laplacian, we impose the following two sets of conditions:
\begin{align}
\lf.s_1 \rt|_{{{x}} \in P_1, {{y}} \in \mathbb{B}_1} &= \lf.\Delta \|{{x}} - {{x}}_{\calI^1_j}\|^m\rt|_{{{x}} \in \mathbb{B}_1},\quad j=1,\hdots,n, \label{eq:interp_constraint}\\
\sum_{j=1}^n \lf.w_j^1({{y}}) \psi_i^1({{x}}) \rt|_{{{x}} \in P_1, {{y}} \in \mathbb{B}_1} &= \lf.\Delta \psi^1_i({{x}})\rt|_{{{x}} \in \mathbb{B}_1},\quad i=1,\hdots,M. \label{eq:poly_constraint}
\end{align}
The interpolant \eqref{eq:rbf_interp} and the two conditions \eqref{eq:interp_constraint}--\eqref{eq:poly_constraint} can be collected into the following block linear system:
\begin{align}
\begin{bmatrix}
A_1 & \Psi_1 \\
\Psi_1^T & O
\end{bmatrix}
\begin{bmatrix}
W_1 \\
\Lambda_1
\end{bmatrix}
=
\begin{bmatrix}
B_{A_1} \\
B_{\Psi_1}
\end{bmatrix},
\label{eq:rbf_linsys}
\end{align}
where
\begin{align}
(A_1)_{ij} &= \|{{x}}_{\calI^1_i} - {{x}}_{\calI^1_j} \|^m,\quad i,j=1,\hdots,n, \\
(\Psi_1)_{ij} &= \psi^1_j({{x}}_{\calI^1_i}),\quad i=1,\hdots,n,\; j=1,\hdots,M,\\
(B_{A_1})_{ij} &= \lf.\Delta \|{{x}} - {{x}}_{\calR^1_j} \|^m \rt|_{{{x}} = {{x}}_{\calI^1_i}},\quad i=1,\hdots,n, \;j=1,\hdots,p_1, \\
(B_{\Psi_1})_{ij} &= \lf.\Delta \psi^1_i({{x}})\rt|_{{{x}} = {{x}}_{\calR^1_j}},\quad i=1,\hdots,M,\; j=1,\hdots,p_1,\\
O_{ij} &= 0,\quad i,j = 1,\hdots,M.
\end{align}
$W_1$ is the \emph{matrix} of overlapped RBF-FD weights, with each column containing the RBF-FD weights for a point ${{y}} \in \mathbb{B}_1$. The linear system \eqref{eq:rbf_linsys} has a unique solution if the nodes in $P_1$ are distinct~\cite{Fasshauer:2007,Wendland:2004}. The matrix of polynomial coefficients $\Lambda_1$ is a set of Lagrange multipliers that enforces the polynomial reproduction constraint \eqref{eq:poly_constraint}, and can be discarded. This procedure is repeated for all points ${{x}} \in X$ that were \emph{not} in the sets $\mathbb{B}$. The weights are then collected into the rows of a sparse matrix $L$, with each stencil populating multiple rows of $L$. For more details, see~\cite{ShankarJCP2017,SFJCP2018}. In this setting, for a given function $f$, the error in approximating $\Delta f$ at a point ${{y}} \in \mathbb{B}$ is given by~\cite{DavydovSchabackMinimal2018}:
\begin{align}
\lf | \Delta f - \sum\limits_{j=1}^n w_j({{y}}) f({{x}}_j) \rt| \leq C({{y}},f) \lf(h({{y}}) \rt)^{\ell-1},
\end{align}
where $h({{y}})$ is the largest distance between the points ${{y}}$ and all points in a given stencil; clearly, the error depends on the polynomial degree $\ell$. Thus, if we wished for a convergence rate of order $\xi$, we would need to set the polynomial degree as $\ell = \xi + 1$. Following~\cite{SFJCP2018}, we set the stencil size $n$ to be $n = 2M + \lf \lfloor \ln(2M) \rt\rfloor$, where $M = {\ell + d \choose d}$. Finally, again from~\cite{SFJCP2018}, we select $\delta$ to be:
\[
 \delta =
  \begin{cases} 
      \hfill 0.7 \hfill & \text{ if $\ell \leq 3$}, \\
      \hfill 0.5 \hfill & \text{ if $4 \leq \ell < 6$}.			
  \end{cases}
\]

\subsection{Numerical solution of the PDE}
\label{sec:pde_disc}

Having discussed the discretization of differential operators using overlapped RBF-FD, we now turn to the discretization of the PDE \eqref{pair} and boundary conditions \eqref{pairbc}. In terms of the solution $v_{0},v_{1}$ to the PDE \eqref{pair}, recall that \eqref{sum} implies that the large time expected neurotransmitter is
\begin{align*}
\lim_{t\to\infty}\E[c(x,t)]
=v(x):=v_{0}(x)+v_{1}(x).
\end{align*}
We denote the numerical approximations of $v_{0}$, $v_{1}$, and $v$ by ${{\mathbb{V}_0}}$, ${{\mathbb{V}_1}}$, and ${{\mathbb{V}}}:={{\mathbb{V}_0}}+{{\mathbb{V}_1}}$.

Now, the problem \eqref{pair}-\eqref{pairbc} is a pair of coupled Helmholtz-type PDEs, with the variables sharing the same Neumann boundary condition on the outer domain boundary, and having different boundary conditions on any inner domain boundaries. Let $X = \{{{x}}_k\}_{k=1}^N$ be some set of nodes in the domain of interest; in our case, we used the node generator from~\cite{SFKSISC2018} to generate these node sets. For convenience, we partition $X$ into two sets:
\begin{align}
X = X_i \cup X_b,
\end{align}
where $X_i$ are the $N_i$ nodes in the interior of the domain, $X_b = [X_b^o \ X_b^i]$ are the $N_b$ boundary nodes; the superscripts $o$ and $i$ indicate outer and inner boundaries, respectively. To solve the PDE, we enforce that the PDEs hold up to and including the domain boundaries, and that the boundary conditions hold true at the boundaries. This allows us to use a set of ghost nodes outside the domain to enforce boundary conditions, giving us additional numerical stability~\cite{SFJCP2018}. Let $X_g = [X_g^o \ X_g^i]$ be the set of ghost nodes associated with the boundary nodes. Thus, in practice, we are operating on an extended node set $\tilde{X}$ given by
\begin{align}
\tilde{X} = X \cup X_g.
\end{align}
Given this partitioning of the node sets, There are many reasonable possible orderings of the unknowns. For simplicity, we choose to order the unknowns as 
\begin{align}
\left({{\mathbb{V}_0}}\right)_{\tilde{X}} &= \begin{bmatrix}
\left({{\mathbb{V}_0}}\right)_i^T & \left({{\mathbb{V}_0}}\right)_b^T & \left({{\mathbb{V}_0}}\right)_g^T
\end{bmatrix}^T, \\
\left({{\mathbb{V}_1}}\right)_{\tilde{X}} &= \begin{bmatrix}
\left({{\mathbb{V}_1}}\right)_i^T & \left({{\mathbb{V}_1}}\right)_b^T & \left({{\mathbb{V}_1}}\right)_g^T
\end{bmatrix}^T.
\end{align}
Note that
\begin{align}
\left({{\mathbb{V}_0}}\right)_b = \begin{bmatrix} ({{\mathbb{V}_0}})_b^o & ({{\mathbb{V}_0}})_b^i \end{bmatrix}^T,
\end{align}
that is, we order ${{\mathbb{V}_0}}$ at the boundary nodes as ${{\mathbb{V}_0}}$ at the outer boundary nodes followed by ${{\mathbb{V}_0}}$ inner boundary nodes; the same ordering is used for ${{\mathbb{V}_1}}$. Of course, for consistency, the same ordering is followed for ghost nodes as well. Now, given an arbitrary function $f$, we know that $L_{\tilde{X}} f_{\tilde{X}} \approx \left.\Delta f\right|_{\tilde{X}}$. Consistent with the above ordering, we partition the sparse matrix $L_{\tilde{X}}$ as follows:
\begin{align}
L_{\tilde{X}} = \begin{bmatrix}
L_{ii} & L_{ib} & L_{ig} \\
L_{bi} & L_{bb} & L_{bg} \\
L_{gi} & L_{gb} & L_{gg}
\end{bmatrix},
\end{align}
where the subscript $ii$ indicates mapping interior nodes to interior nodes, and so forth.  We denote the discretized Neumann operator corresponding to the outer boundary by $D^o$ (a sparse matrix). We also require a matrix for Neumann boundary conditions on the inner boundaries; denote this by $D^i$. The Dirichlet boundary conditions on the interior do not require a special boundary operator, and can be enforced using identity matrices and zero matrices. We partition these matrices similar to $L$, with the exception that there are no rows mapping any nodes to interior points. To aid in the discretization of the PDE, define the boundary operators for ${{\mathbb{V}_0}}$ to be
\begin{align}
D^0_{bi} &=
\begin{bmatrix}
(D^o_{bi})^T & (0^i_{bi})^T
\end{bmatrix}^T, \\
D^0_{bb} &=
\begin{bmatrix}
(D^o_{bb})^T & (I^i_{bb})^T
\end{bmatrix}^T, \\
D^0_{bg} &= 
\begin{bmatrix}
(D^o_{bg})^T & (0^i_{bg})^T
\end{bmatrix}^T,
\end{align}
which corresponds to enforcing a Neumann boundary condition on the outer boundary, and a Dirichlet condition on the inner boundary. Similarly, the boundary operators for ${{\mathbb{V}_1}}$ are:
\begin{align}
D^1_{bi} &=
\begin{bmatrix}
(D^o_{bi})^T & (D^i_{bi})^T
\end{bmatrix}^T, \\
D^1_{bb} &=
\begin{bmatrix}
(D^o_{bb})^T & (D^i_{bb})^T
\end{bmatrix}^T, \\
D^1_{bg} &= 
\begin{bmatrix}
(D^o_{bg})^T & (D^i_{bg})^T
\end{bmatrix}^T,
\end{align}
which corresponds to Neumann boundary conditions on both outer and inner boundaries. With this notation in place, the fully discretized (block) system of PDEs and boundary conditions is given by:
{\small
\begin{align}
\underbrace{
\begin{bmatrix}
L_{ii} - \beta I_{ii} & L_{ib} & L_{ig} & \alpha I_{ii} & 0 & 0 \\
L_{bi} & L_{bb} - \beta I_{bb} & L_{bg} & 0 & \alpha I_{bb} & 0 \\
D^0_{bi} & D^0_{bb} & D^0_{bg} & 0 & 0 & 0 \\
\beta I_{ii} & 0 & 0 & L_{ii}-\alpha I_{ii} & L_{ib} & L_{ig} \\
0 & \beta I_{bb} & 0 & L_{bi} & L_{bb} - \alpha I_{bb} & L_{bg} \\
0 & 0 & 0 & D^1_{bi} & D^1_{bb} & D^1_{bg}
\end{bmatrix}}_{\tilde{L}}
\underbrace{\begin{bmatrix}
\left({{\mathbb{V}_0}}\right)_i \\
\left({{\mathbb{V}_0}}\right)_b \\
\left({{\mathbb{V}_0}}\right)_g \\
\left({{\mathbb{V}_1}}\right)_i \\
\left({{\mathbb{V}_1}}\right)_b \\
\left({{\mathbb{V}_1}}\right)_g \\
\end{bmatrix}}_{\tilde{V}}
=
\underbrace{
\begin{bmatrix}
0 \\
0 \\
0 \\
0 \\
0 \\
\tilde{\Theta}
\end{bmatrix}}_{\tilde{R}},
\end{align}}
where $\tilde{\Theta} = [0 \ \Theta]^T$, and $\Theta$ is a vector with length equal to the number of inner boundary points with each entry as the constant $\frac{\beta}{\alpha + \beta}$. In the above system, $I$ refers to the identity matrix, with the subscripts indicating the size. $0$ in each case refers to a matrix of zeros; since the size of these zero blocks can be inferred from the rows in which they sit, we leave out subscripts for brevity. The block matrix $\tilde{L}$ is very sparse, since its individual blocks are themselves sparse matrices. Note that this matrix has dimensions $2(N + N_b)$. Since the continuous problem is well-posed, a good discretization of the continuous problem should result in a matrix $\tilde{L}$ that does not have a nullspace. In practice, this translates into a requirement on both the node sets used to discretize the domain, and on the RBF-FD discretizations of the Laplacian and Neumann operators. We are able to invert the matrix $\tilde{L}$ using Matlab's built-in sparse direct solver without any difficulties. We choose the sparse direct solver for simplicity, rather than using a preconditioned iterative solve. We leave such extensions for future work.

\subsection{Numerical Results}
\label{sec:num_results}

In this section, we present two types of numerical results. First, we show results that validate our numerical method on the spherically symmetric problem of section~\ref{shell} for which we have the analytical solution \eqref{explshell}. Then, using this validation test to guide our choices, we present a numerical investigation of our asymptotic results in two different domain geometries.

\subsubsection{Validation of numerical methods}
\label{sec:num_valid}

We first consider the spherically symmetric problem of section~\ref{shell} in which the domain is the unit sphere with a spherical hole of radius $\eps>0$ centered at the origin. In this scenario, the large-time expected neurotransmitter $v=v_{0}+v_{1}$ is exactly constant in space and this constant is given explicitly in \eqref{explshell}. To validate our numerical methods, we thus want to show that (i) the numerical solution is approximately constant in space and (ii) this approximately constant value of the numerical solution converges to the value of \eqref{explshell}.

Before checking these two points, we note that since this is a problem in three space dimensions and spatial refinement is relatively expensive, we instead measure convergence by order refinement. In the overlapped RBF-FD method, this is accomplished by increasing $\xi$; we use $\xi \in [2,5]$ for this study. In this case, we expect the error to decrease as $\xi$ increases, regardless of the hole radius $\eps$. However, when we first attempted this study for holes of decreasing radius, we observed that the errors were higher for smaller holes, primarily due to our use of a fixed $h$ in the spatial discretization. After some experimentation, we found that maintaining a fixed ratio of $h/\eps$ was critical to ensuring that the spatial errors were independent of hole radius (we take $h/\eps \approx 1/7$). We continue to use this refinement strategy in section~\ref{sec:num_asymp}.

\begin{figure}
	\centering
		\includegraphics[width=12cm]{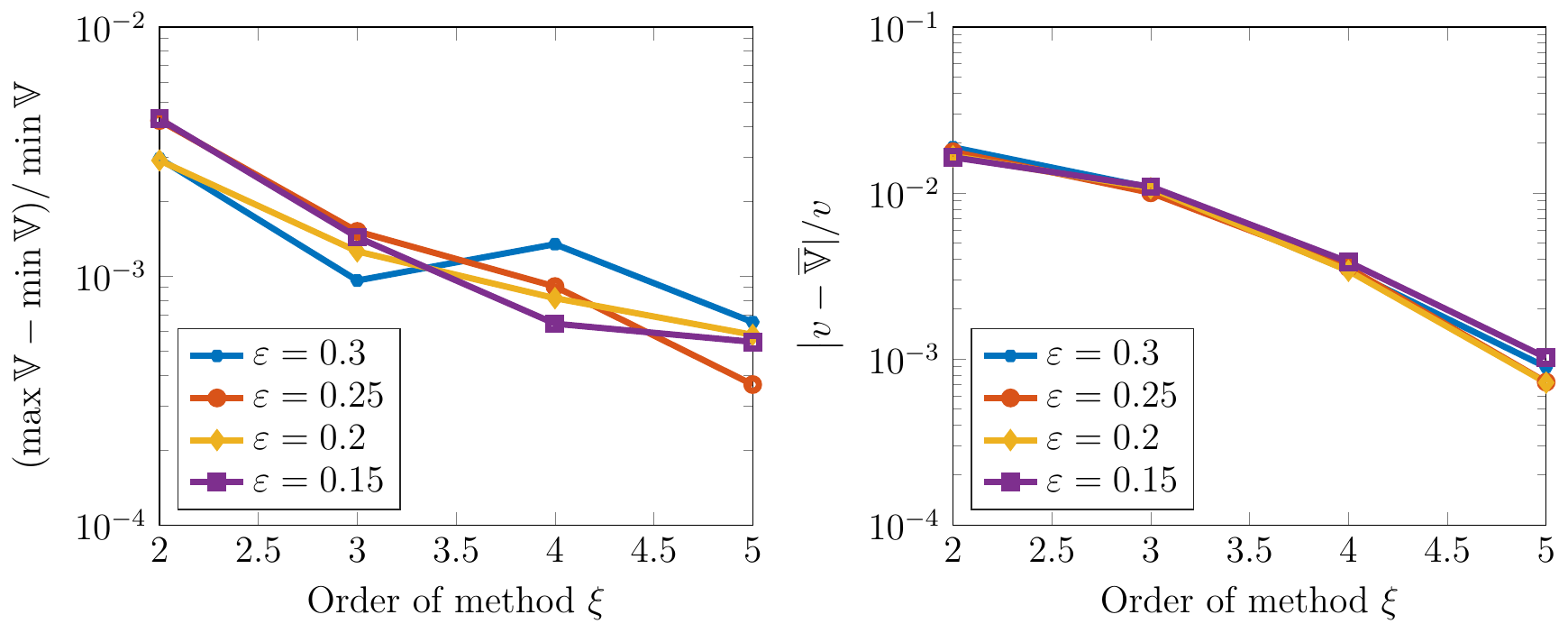}
		\caption{For the spherically symmetric example in section~\ref{shell}, the left panel confirms that the numerical solution is approximately spatially constant. The right panel confirms that this approximately constant numerical solution approaches the exactly constant analytical solution given in \eqref{explshell}. 
		See the text for details.}
	\label{fig:conv_test}
\end{figure}

To show that the numerical solution is approximately constant in space, the left panel of Figure~\ref{fig:conv_test} plots the following measure of spatial variation,
\begin{align*}
\frac{\max\mathbb{V}-\min\mathbb{V}}{\min\mathbb{V}},
\end{align*}
where $\mathbb{V}=\mathbb{V}_{0}+\mathbb{V}_{1}$ is the numerical approximation of $v=v_{0}+v_{1}$ and the maximums and minimums are taken over all interior nodes. This plot confirms that the deviation from a spatially constant solution decays under order refinement. Indeed,  this plot shows that the maximum and minimum of the solution are well within $1\%$ for the lowest order method and well within $0.1\%$ for the highest order method, and this holds regardless of hole radius $\eps$. 

To show that the (approximately constant) value of the numerical solution approaches the analytical value $v=v_{0}+v_{1}$ under order refinement, the right panel of Figure~\ref{fig:conv_test} plots the normalized error, $|v-\overline{\mathbb{V}}|/v$, where $\overline{\mathbb{V}}$ denotes the value of the numerical solution averaged over all interior grid points. This plot shows that the error drops off under order refinement, reaching much less than $0.1\%$ for the highest order method.

\subsubsection{Investigation of asymptotics}
\label{sec:num_asymp} 

Having validated our numerical methods in section~\ref{sec:num_valid}, we employ our numerical discretization to verify some of our asymptotic results. Specifically, we verify the asymptotic result in \eqref{efo} in section~\ref{sync} in the following two geometries. The first is for a spherical domain of radius $R_{0}=1$ with two holes of radius $\eps$ centered at the Cartesian coordinates $(-0.3,-0.3,-0.3)\in\R^{3}$ and $(0.3,0.3,0.3)\in\R^{3}$. The second is a spherical domain of radius $R_{0}=1.2$ with the same holes. To ensure that the discretization error is small, we use $\xi = 4$ (though similar results were obtained with $\xi=3$ and $\xi=5$). 

\begin{figure}
	\centering
		\includegraphics[width=12cm]{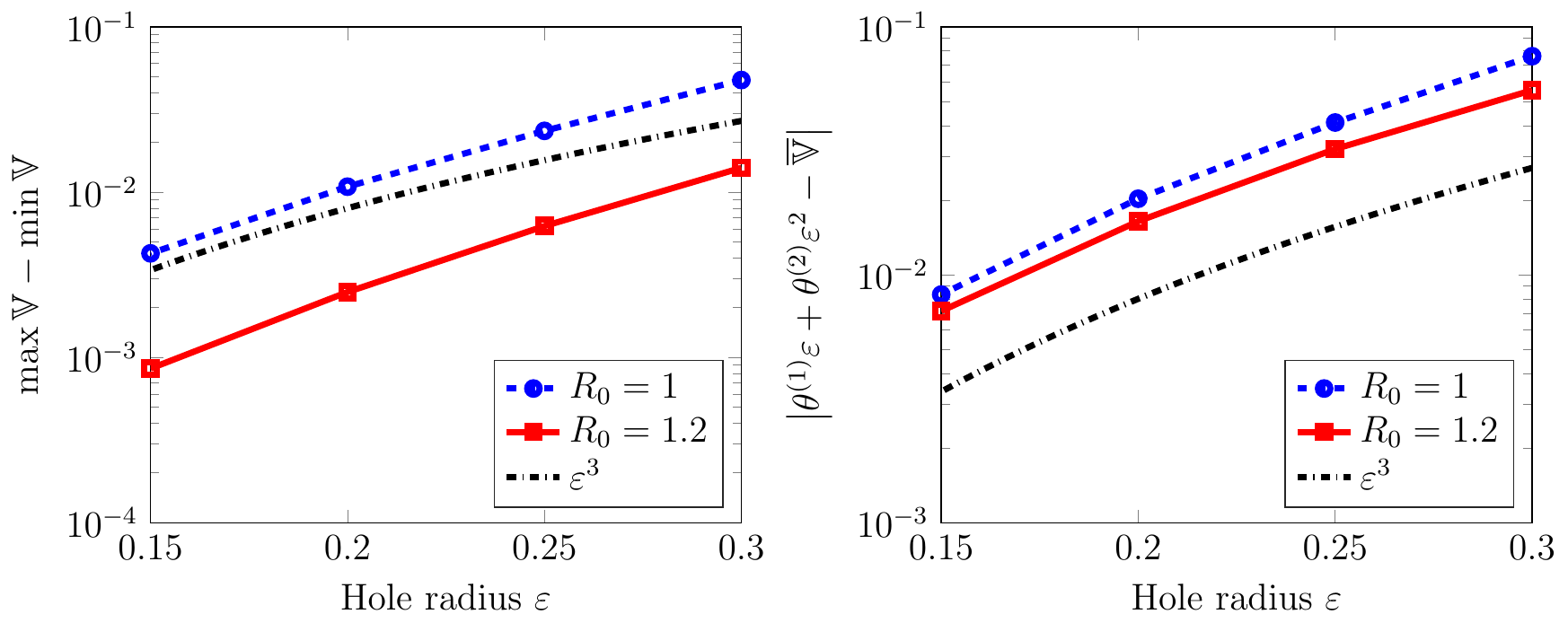}
		\caption{
		For the synchronous firing example of section~\ref{sync} with two interior holes, the left panel confirms that the numerical solution is approximately spatially constant, as predicted by the asymptotics. The right panel shows that the difference between the numerical solution and the two-term asymptotic expansion is $\O(\eps^{3})$, as predicted by the asymptotics. The blue dashed and red solid curves correspond to a spherical domain of radius $R_{0}=1$ and $R_{0}=1.2$, respectively. See the text for details.}
	\label{fig:asymp_test}
\end{figure}

Now, the asymptotic result in \eqref{efo} predicts that the solution is approximately constant outside an $\O(\eps)$ neighborhood of each of the holes. Letting ${{\mathbb{V}}}={{\mathbb{V}_0}} + {{\mathbb{V}_1}}$ denote the numerical estimate to $v=v_{0}+v_{1}$ where $v_{0},v_{1}$ are the (unknown) exact solutions to \eqref{pair}-\eqref{pairbc}, the left panel of Figure~\ref{fig:asymp_test} confirms that the numerical solution is approximately spatially constant and that the deviation from constant decays as $\eps$ decays. Indeed, the left panel plots the difference, $\max\mathbb{V}-\min\mathbb{V}$, where the maximum and minimum are taken over all interior grid points, and this difference is approximately $\O(\eps^{3})$ or $\O(\eps^{4})$, which agrees with the asymptotic result in \eqref{efo}.

The right panel of Figure~\ref{fig:asymp_test} plots the absolute value of the difference between the asymptotic estimate, $\eps\theta^{(1)}+\eps^{2}\theta^{(2)}$, in \eqref{efo} and the numerical solution averaged over all interior grid points, $\overline{\mathbb{V}}$. In agreement with \eqref{efo}, this difference is $\O(\eps^{3})$. Indeed, the lines of best fit to the $R_{0}=1$ and $R_{0}=1.2$ curves in the semi-log plot of the right panel yields respective slopes of 3.19 and 2.96.

\section{Discussion}

We have used asymptotic analysis to study the large-time mean of a stochastic PDE that models the concentration of a neurotransmitter in the extracellular space of a brain region. This analysis yielded the first two terms in an asymptotic expansion of the concentration, which shows the first and second order effects of the many details and parameters in the problem. We developed a numerical solution method to verify our asymptotic results in some special cases.

Our asymptotic analysis assumed that the nerve varicosities are much smaller than the typical distance between them. Following \cite{lawley18prob}, the distance between varicosities varies, but for serotonin there are approximately $2.6\times10^{6}$ varicosities per cubic millimeter \cite{soghomonian87}, which implies a distance between varicosities of approximately $l_{1}=7\,\mu\text{m}$. In Figure 1 of \cite{pasik82}, varicosities are separated by about $l_{1}=20\,\mu\text{m}$ on average. The average radius of a varicosity in \cite{soghomonian87} is approximately $l_{0}=.3\,\mu\text{m}$. Hence, we estimate that the ratio $\eps=l_{0}/l_{1}$ is indeed small,
\begin{align*}
.015=.3/20\le\eps\le .3/7 \approx.04.
\end{align*}

Of course, our model ignores various features of the actual biological system. For one, the rate of neurotransmitter release is likely not constant when a neuron is firing, and the varicosities do not perfectly absorb neurotransmitter when the neuron is not firing. In fact, a more detailed model might allow the neurotransmitter release rate and/or firing rate to depend on the extracellular concentration, and these may be interesting avenues for further research. Furthermore, we have ignored the presence of obstacles in the extracellular space which restrict the diffusive path of neurotransmitter (such as the actual neurons in the extracellular space; we have merely included the bulbous varicosities on a neuron). Indeed, the extracellular space comprises about 20\% of brain volume and is a very tortuous and highly connected foamlike structure formed from the interstices of cell surfaces \cite{sykova2008}. However, using the diffusion equation to model biochemical concentrations in the extracellular space in the brain has a long history \cite{nicholson1981} and its validity is well-established for length scales of approximately $2$ to $200\,\mu\text{m}$ \cite{nicholson2001, nicholson2017, sykova2008}. The tortuosity of the extracellular space is often modeled by a reduction in the effective diffusion coefficient compared to diffusion in free solution \cite{nicholson2017}. The tortuosity could also potentially be modeled by a space-dependent diffusion coefficient.

Perhaps more importantly, our model does not include non-neuronal cells which can take up neurotransmitter. While the brain contains $10^{11}$ neurons, it actually contains at least ten times as many glial cells \cite{feldman}. Furthermore, some glial cells are known to take up neurotransmitter \cite{daws13,daws05}, though their neurotransmitter absorption rate is likely not as fast as that of nerve varicosities. In fact, a recent study combining experiments and mathematical modeling suggests that this uptake pathway is important for serotonin \cite{wood14}. An important extension of our model and analysis would be to understand how these glia affect the extracellular neurotransmitter concentration. We expect that the steady-state neurotransmitter concentration would no longer be spatially constant if we included this uptake by glial cells. However, it remains to determine the magnitude of the spatial variation that would be introduced by this uptake.

An additional avenue for future work is investigating the variance of the neurotransmitter concentration in the extracellular space. Previous analysis of a simplified model in one space dimension showed that this variance is approximately constant outside a small neighborhood of a varicosity, but quantitative estimates are lacking for models in three space dimensions. This analysis would require significant extensions of the asymptotic and numerical methods employed in the present work, since calculating the variance of a three-dimensional stochastic model requires solving a certain boundary value problem in six space dimensions which couples at the boundaries to the mean of the stochastic model \cite{lawley16bvp}.

\medskip
{\bf Acknowledgements.}  SDL thanks Janet Best, Fred Nijhout, and Michael Reed for stimulating discussions. 

\bibliography{library}
\bibliographystyle{unsrt}
\end{document}